\documentclass[12]{article}

\usepackage{fleqn}
\usepackage{graphicx,epsf,psfrag}
\usepackage{amsthm,amsfonts,amssymb,amsmath}

\begin{document}
\Large
\begin{center}
{\bf Cayley-Dickson Algebras and Finite Geometry}
\end{center}
\vspace*{.1cm}
\large
\begin{center}
Metod Saniga,$^{1,2}$ Fr\'ed\'eric Holweck$^{3}$ and Petr Pracna$^{4}$ 
\end{center}
\vspace*{-.4cm}
\normalsize
\begin{center}
$^{1}$Institute for Discrete Mathematics and Geometry,
Vienna University of Technology\\ Wiedner Hauptstra\ss e 8--10,
A-1040 Vienna, Austria\\
(metod.saniga@tuwien.ac.at)

\vspace*{.2cm}

$^{2}$Astronomical Institute, Slovak Academy of Sciences\\
SK-05960 Tatransk\' a Lomnica, Slovak Republic\\
(msaniga@astro.sk) 

\vspace*{.2cm}

$^{3}$Laboratoire IRTES/M3M, Universit\'e de Technologie de Belfort-Montb\'eliard\\ 
F-90010 Belfort, France\\ (frederic.holweck@utbm.fr) 

\vspace*{.2cm}

$^{4}$J. Heyrovsk\' y Institute of Physical Chemistry, v.\,v.\,i.\\  Academy of Sciences of the Czech Republic \\ Dolej\v skova 3, CZ-18223 Prague, Czech Republic\\
(pracna@jh-inst.cas.cz)

\end{center}

\vspace*{-.3cm} \noindent \hrulefill

\vspace*{-.0cm} \noindent {\bf Abstract}

\noindent Given a $2^N$-dimensional Cayley-Dickson algebra, where $3 \leq N \leq 6$, we first observe  that the multiplication table of its imaginary units $e_a$, $1 \leq a \leq 2^N -1$, is encoded in the properties of the projective space PG$(N-1,2)$ if one regards these imaginary units as points and distinguished triads of them $\{e_a, e_b, e_c\}$, $1 \leq a < b <c \leq 2^N -1$ and $e_ae_b = \pm e_c$, as lines. This projective space is seen to feature two distinct kinds of lines according as $a+b = c$ or  $a+b \neq c$. Consequently, it also exhibits (at least two) different types of points in dependence on how many lines of either kind pass through each of them. In order to account for such partition of the PG$(N-1,2)$, the concept of Veldkamp space of a finite point-line incidence structure is employed. The corresponding point-line incidence structure is found to be a binomial $\left({N+1 \choose 2}_{N-1}, {N+1 \choose 3}_{3}\right)$-configuration ${\cal C}_N$; in particular, ${\cal C}_3$ (octonions)  is isomorphic to the Pasch $(6_2,4_3)$-configuration, ${\cal C}_4$ (sedenions) is the famous Desargues $(10_3)$-configuration, ${\cal C}_5$ (32-nions) coincides with the Cayley-Salmon $(15_4,20_3)$-configuration found in the well-known Pascal mystic hexagram and ${\cal C}_6$ (64-nions) is identical with a particular
$(21_5,35_3)$-configuration that can be viewed as four triangles in perspective from a line where the points of perspectivity of six pairs of
them form a Pasch configuration. We also draw attention to a remarkable nesting pattern formed by these configurations, where ${\cal C}_{N-1}$ occurs as a geometric hyperplane of ${\cal C}_N$. Finally, a brief examination of the structure of generic ${\cal C}_N$ leads to a conjecture that
${\cal C}_N$ is isomorphic to a combinatorial Grassmannian of type $G_2(N+1)$.

\vspace*{.3cm}


\vspace*{-.2cm} \noindent \hrulefill


\section{Introduction} 

As it is well known (see, e.\,g., \cite{jac,scha}), the Cayley-Dickson algebras represent a nested sequence $A_0, A_1, A_2,\ldots, A_N, \ldots$ of $2^N$-dimensional (in general non-associative) $\mathbb{R}$-algebras with $A_N \subset A_{N+1}$, where $A_0 = \mathbb{R}$ and where for any $N \geq 0$ $A_{N+1}$ comprises  all ordered pairs of elements from $A_{N}$ with conjugation defined by
\begin{equation}
(x,y)^{\ast} = (x^{\ast}, -y) 
\label{eq1}
\end{equation}
and multiplication usually by
\begin{equation}
(x,y)(X,Y) = (xX - Yy^{\ast}, x^{\ast}Y + Xy). 
\label{eq2}
\end{equation}        
Every finite-dimensional algebra (see, e.\,g., \cite{scha}) is basically defined by the multiplication rule of its basis. 
The basis elements (or units) $e_0, e_1, e_2, \ldots, e_{2^{N+1} -1}$ of $A_{N+1}$, $e_0$ being the real basis element (identity), can be chosen in various ways. Our preference is the {\it canonical} basis
\begin{align*}
e_0     &= (e_0,0), & e_1       &= (e_1,0), & e_2       &= (e_2,0), && \ldots, &&   e_{2^N - 1}    &= (e_{2^N - 1},0), &~~~~~& \\
e_{2^N} &= (0,e_0), & e_{2^N+1} &= (0,e_1), & e_{2^N+2} &= (0,e_2), && \ldots, && e_{2^{N+1} - 1} &= (0,e_{2^{N} - 1}),&~~~~~~&
\end{align*}
where, by abuse of notation (see, e.\,g., \cite{mor}), the same symbols are also used for the basis elements of $A_N$.
This is because the paper essentially focuses on {\it multiplication} properties of basis elements and the canonical basis seems to
display most naturally the inherent symmetry of this operation. For,
in addition to revealing the nature of the Cayley-Dickson recursive process, it also implies that for both $a$ and $b$ being non-zero we have
$e_a e_b = \pm e_{a \oplus b},$  where the symbol `$\oplus$' denotes `{\it exclusive or}' of the binary representations of $a$ and $b$ (see, e.\,g., \cite{arndt}). 
From the above expressions and Eqs.~(\ref{eq1}) and (\ref{eq2}) one can readily find the product of any two distinct units of $A_{N+1}$
if the multiplication properties of those of $A_{N}$ are given.
Such products are usually expressed/presented in a tabular form, and we shall also follow this tradition here. All the multiplication tables we made use of were computed for us by J\"org Arndt; the computer program is named {\tt cayley-dickson-demo.cc} and is freely available from his web-site {\tt http://jjj.de/fxt/demo/arith/index.html}.

Employing such a multiplication table of $A_N$, $N \geq 2$, it can be verified that the $2^N-1$ imaginaries $e_a$, $1 \leq a \leq 2^N -1$, form ${2^N-1 \choose 2}/3$ distinguished sets each of which comprises three different units $\{e_a, e_b, e_c\}$ that satisfy equation 
\begin{equation}
e_a e_b = \pm e_c,
\label{eq4}
\end{equation}
and where each unit is found to belong to $2^{N-1} - 1$ such sets.\footnote{Although these properties will explicitly be illustrated only for the cases $3 \leq N \leq 6$, due to the nested structure of the algebras they must be exhibited by any other higher-dimensional $A_N$.} 
Regarding the imaginaries as points and their distinguished triples as lines, one gets a point-line incidence geometry where every line has three points and through each point there pass $2^{N-1} - 1$ lines and which is isomorphic to PG$(N-1,2)$, the $(N-1)$-dimensional projective space over the smallest Galois field $GF(2)$ (see also \cite{baez} for $N=3$ and \cite{culb} for $N=4$).
Let us assume, without loss of generality, that the elements in any distinguished triple $\{e_a, e_b, e_c\}$ of $A_N$ are ordered in such a way that $a < b < c$.  Then, for $N \geq 3$, we can naturally speak about two different kinds of triples and, hence, two distinct kinds of lines of the associated $2^N$-nionic PG$(N-1,2)$, according as $a+b = c$ or  $a+b \neq c$; in what follows a line of the former/latter kind will be called ordinary/defective. This stratification of the line-set of the PG$(N-1,2)$ induces a similar partition of the point-set of the latter space into several types, where a point of a given type is characterized by the same number of lines of either kind that pass through it.  

Obviously, if our projective space PG$(N-1,2)$ is regarded as an {\it abstract} geometry {\it per se}, every point and/or every line in it has the same footing.  
So, to account for the above-described `refinement' of the structure of our $2^N$-nionic PG$(N-1,2)$, it turns out to be necessary to find a representation of this space where each point/line is ascribed a certain `internal' structure, which at first sight may seem to be quite a challenging task. 
To tackle this task successfully, we need to introduce a few notions/concepts from the realm of finite geometry. 

We start with a finite {\it point-line incidence structure} $\mathcal{C} = (\mathcal{P},\mathcal{L},I)$ where $\mathcal{P}$ and $\mathcal{L}$ are, respectively, finite sets of points and lines and where incidence $I \subseteq \mathcal{P} \times \mathcal{L}$ is a binary relation indicating which point-line pairs are incident (see, e.\,g., \cite{shult}).    
Here, we shall only be concerned with specific point-line incidence structures called {\it configurations} \cite{grun}. A $(v_r,b_k)$-configuration is a $\mathcal{C}$ where: 1) $v = \vert \mathcal{P} \vert$ and $b = \vert \mathcal{L}\vert$, 2) every line has $k$ points and every point is on $r$ lines, and 3) two distinct lines intersect in at most one point and every two distinct points are joined by at most one line; a configuration where $v=b$ and $r=k$ is called symmetric (or balanced), and usually denoted as a $(v_r)$-configuration.
A $(v_r,b_k)$-configuration with $v = {r+k-1 \choose r}$ and $b = {r+k-1 \choose k}$ is called a $binomial$ configuration.
Next it comes a {\it geometric hyperplane} of $\mathcal{C} = (\mathcal{P},\mathcal{L},I)$, which is a proper subset of $\mathcal{P}$ such that a line from $\mathcal{C}$ either lies fully in the subset, or shares with it only one point. 
 If $\mathcal{C}$ possesses geometric hyperplanes, then one can define the {\it Veldkamp space} of $\mathcal{C}$,  $\mathcal{V}(\mathcal{C})$, as follows \cite{buec}: (i) a point of $\mathcal{V}(\mathcal{C})$ is a geometric hyperplane of  $\mathcal{C}$
and (ii) a line of $\mathcal{V}(\mathcal{C})$ is the collection $H'H''$ of all geometric hyperplanes $H$ of $\mathcal{C}$  such that $H' \cap H'' = H' \cap H = H'' \cap H$ or $H = H', H''$, where $H'$ and $H''$ are distinct geometric hyperplanes. If each line of $\mathcal{C}$ has three points and $\mathcal{C}$ `behaves well,' a line of $\mathcal{V}(\mathcal{C})$ is also of size three and can equivalently be defined as $\{H', H'', \overline{H' \Delta H''}\}$, where the symbol $\Delta$ stands for the symmetric difference of the two geometric hyperplanes and an overbar denotes the complement of the object indicated.
From its definition it is obvious that  $\mathcal{V}(\mathcal{C})$ is well suitable for our needs because its points, being themselves {\it sets} of points, have different `internal' structure and so, in general, they can no longer be on the same par; clearly, the same applies to the lines $\mathcal{V}(\mathcal{C})$. Our task thus basically boils down to finding such $\mathcal{C}_N$ whose  $\mathcal{V}(\mathcal{C}_N)$ is isomorphic to PG$(N-1,2)$ and completely reproduces its $2^N$-nionic fine structure. This will be carried out in great detail for the first four non-trivial cases, $3 \leq N \leq 6$, which, when combined with the two trivial cases ($N=1,2$), will provide us with sufficient amount of information to guess a general pattern.  

The paper is organized as follows. In Sect.~2 it is shown that ${\cal C}_3$ (octonions)  is isomorphic to the Pasch $(6_2,4_3)$-configuration, which plays a key role in classifying Steiner triple systems. In Sect.~3 one demonstrates that ${\cal C}_4$ (sedenions) is nothing but the famous Desargues $(10_3)$-configuration. In Sect.~4 our ${\cal C}_5$ (32-nions) is shown to be identical with the Cayley-Salmon $(15_4,20_3)$-configuration found in the well-known Pascal mystic hexagram.  In Sect.~5 we find that ${\cal C}_6$ corresponds to a particular $(21_5,35_3)$-configuration encompassing seven distinct copies of the Cayley-Salmon $(15_4,20_3)$-configuration as geometric hyperplanes. In Sect.~6 some rudimentary properties of the generic ${\cal C}_N \cong \left({N+1 \choose 2}_{N-1}, {N+1 \choose 3}_{3}\right)$-configuration are outlined and its isomorphism to a combinatorial Grassmannian of type $G_2(N+1)$ is conjectured. Finally, Sect.~7 is reserved for concluding remarks.

\section{Octonions and the Pasch $(6_2,4_3)$-configuration}
From the nesting property of the Cayley-Dickson construction of $A_N$ it is obvious that the smallest non-trivial case to be addressed is
$A_3$, the algebra of octonions, whose multiplication table is presented in Table 1.

\begin{table}[h]
\begin{center}
\caption{The multiplication table of the imaginary unit octonions $e_a$, $1 \leq a \leq 7$. For the sake of simplicity, in what follows we shall employ a short-hand notation $e_a \equiv a$; likewise for the real unit $e_0 \equiv 0$. There are also delineated multiplication tables corresponding to the distinguished nested sequence of sub-algebras of complex numbers ($a=1$, the upper left corner) and quaternions ($1 \leq a \leq 3$, the upper left $3 \times 3$ square).} 
\vspace*{0.4cm}
\begin{tabular}{||r|r|rr|rrrr||}
\hline \hline
$*$   &  1    &  2    &  3    &  4    &  5    &  6    &  7   \\
\hline
   1  &  $-$0 &  $-$3 &  +2   &  $-$5 &  +4   &  +7   &  $-$6 \\
\hline
   2  &  +3   &  $-$0 &  $-$1 &  $-$6 &  $-$7 &  +4   &  +5   \\ 
   3  &  $-$2 &  +1   &  $-$0 &  $-$7 &  +6   &  $-$5 &  +4   \\
\hline
   4  &  +5   &  +6   &  +7   &  $-$0 &  $-$1 &  $-$2 &  $-$3 \\
   5  &  $-$4 &  +7   &  $-$6 &  +1   &  $-$0 &  +3   &  $-$2 \\
   6  &  $-$7 &  $-$4 &  +5   &  +2   &  $-$3 &  $-$0 &  +1   \\
   7  &  +6   &  $-$5 &  $-$4 &  +3   &  +2   &  $-$1 &  $-$0 \\
 \hline \hline
\end{tabular}
\end{center}
\end{table}

\noindent
The above-given multiplication table implies the existence of the following seven {\it distinguished} trios of imaginary units:

$$
\{1,2,3\}, \{1,4,5\}, \{1,6,7\}, 
$$
$$
\{2,4,6\}, \{2,5,7\}, 
 $$
$$
\{3,4,7\}, \{3,5,6\}.
$$

\begin{figure}[t]
	\centering
	\includegraphics[width=4.0truecm]{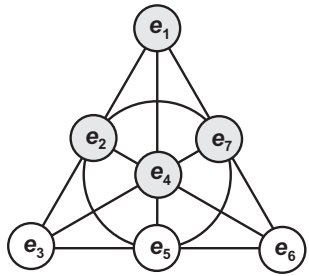}
	\caption{An illustration of the structure of PG$(2,2)$, the Fano plane, that provides the multiplication law for octonions (see, e.\,g., \cite{baez}). The points of the plane are seven small circles. The lines are represented by triples of circles located on the sides of the triangle, on its altitudes, and by the triple lying on the big circle. The three imaginaries lying on the same line satisfy Eq.\,(\ref{eq4}).}
	\label{f1}
\end{figure}

\noindent
Regarding the seven imaginary units as points and the seven distinguished triples of them as lines, we obtain a point-line incidence structure where each line has three points and, dually, each point is on three lines, and which is isomorphic to the smallest projective plane PG$(2,2)$, often called the Fano plane, depicted in Figure \ref{f1}.
It is then readily seen that we have six ordinary lines, namely
$$
\{1,2,3\}, \{1,4,5\}, \{1,6,7\},
$$
$$
\{2,4,6\}, \{2,5,7\}, 
 $$
$$
\{3,4,7\}, 
$$
and only single defective one, viz.
$$
\{3,5,6\}.
$$
Similarly, our octonionic PG$(2,2)$ features two distinct types of points. 
A type-one point is such that two lines passing through it are ordinary, the remaining one being defective; such a point lies in the set 
$$\{3,5,6\} \equiv \alpha.$$
A type-two point is such that every line passing through it is ordinary; such a point belongs  to the set 
$$\{1,2,4,7\} \equiv \beta,$$  
which is highlighted by gray color in Figure \ref{f1}.

\begin{figure}[pth!]
	\centering
	\includegraphics[width=2.5truecm]{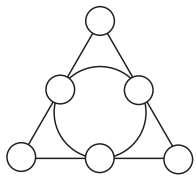}
	\caption{An illustrative portrayal of the Pasch configuration: circles stand for  its points, whereas its lines are represented by triples of points on common straight segments (three) and the triple lying on a big circle.}
\label{f2}
\end{figure}

A configuration $\mathcal{C}_3$ whose Veldkamp space reproduces the above-described partitions of points and lines of PG$(2,2)$ is, as we will soon see, nothing but the well-known {\it Pasch} $(6_2,4_3)$-configuration, $\mathcal{P}$. This configuration, which plays a very important role in classifying Steiner triple systems (see, e.\,g., \cite{enc}), is depicted in Figure \ref{f2} in a form showing an automorphism of order three; it also lives in the Fano plane and, as it is readily seen by comparing Figures \ref{f1} and \ref{f2}, it can be obtained from the latter by removal of any of its seven points and all the three lines passing through it.

\begin{figure}[pth!]
	\centering
	\includegraphics[width=5truecm]{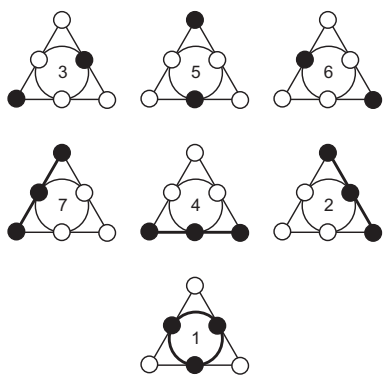}
	\caption{The seven geometric hyperplanes of the Pasch configuration. The hyperplanes are labelled by imaginary units of octonions in such a way that --- as it is obvious from the next figure --- the seven lines of the Veldkamp space of the Pasch configuration are identical with the seven distinguished triples of units, that is with the seven lines of the PG$(2,2)$ shown in Figure \ref{f1}.}
\label{f3}	
\end{figure}

\begin{figure}[pth!]
	\centering
	\includegraphics[width=7.5truecm]{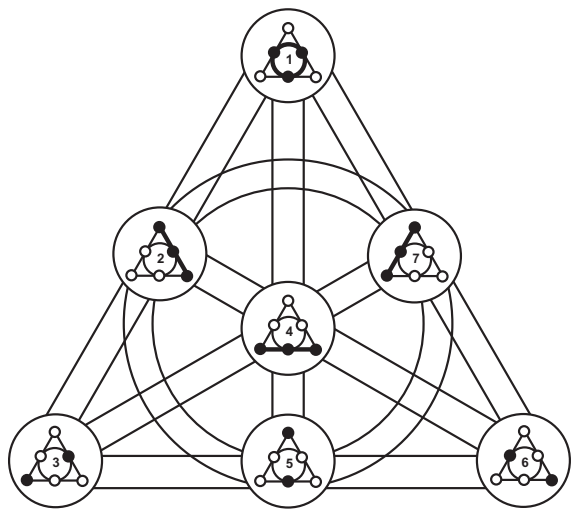}
	\caption{A unified view of the seven Veldkamp lines of the Pasch configuration. The reader can readily verify that for any three geometric hyperplanes lying on a given line of the Fano plane, one is the complement of the symmetric difference of the other two.}
\label{f4}	
\end{figure}

In order to see that $\mathcal{V}(\mathcal{P}) \cong$ PG$(2,2)$ we shall first show, using our diagrammatical representation of $\mathcal{P}$, all seven geometric hyperplanes of $\mathcal{P}$ --- Figure \ref{f3}. We see that they are indeed of two different forms, of cardinality three and four. A member of the former set comprises two points at maximum distance from each other. Such geometric hyperplane corresponds to a type-one (or $\alpha$-) point of PG$(2,2)$. A member of the latter set features three points on a common line;  such a geometric hyperplane of $\mathcal{P}$ corresponds to a type-two (or $\beta$-) point of our PG$(2,2)$.
The seven lines of $\mathcal{V}(\mathcal{P})$ are illustrated in a compact diagrammatic form in Figure \ref{f4}; as it is easily discernible, each of the six ordinary lines is of the form $\{\alpha, \beta, \beta\}$, whilst the remaining defective one has the $\{\alpha, \alpha, 
\alpha\}$ shape.

\section{Sedenions and the Desargues $(10_3)$-configuration}
Our next focus is on $A_4$, the sedenions, whose basic multiplication properties are summarized in Table 2.

\begin{table}[h]
\centering
\caption{The multiplication table of the imaginary unit sedenions $e_a$, $1 \leq a \leq 15$. As in the previous section, we shall employ a short-hand notation $e_a \equiv a$; likewise for the real unit $e_0 \equiv 0$. There are also shown multiplication tables corresponding to the distinguished nested sequence of sub-algebras starting with complex numbers ($a=1$), quaternions ($1 \leq a \leq 3$) and octonions ($1 \leq a \leq 7$).} 
\vspace*{0.4cm}
\resizebox{\columnwidth}{!}{%
\begin{tabular}{||r|r|rr|rrrr|rrrrrrrr||}
\hline \hline
   $*$   &  1    &  2    &  3    &  4    &  5    &  6    &  7    &  8    &  9    & 10    & 11    & 12    & 13    & 14    & 15    \\
\hline
   1  &  $-$0 &  $-$3 &  +2   &  $-$5 &  +4   &  +7   &  $-$6 &  $-$9 &  +8   & +11   & $-$10 & +13   & $-$12 & $-$15 & +14   \\
\hline
   2  &  +3   &  $-$0 &  $-$1 &  $-$6 &  $-$7 &  +4   &  +5   & $-$10 & $-$11 &  +8   &  +9   & +14   & +15   & $-$12 & $-$13 \\
   3  &  $-$2 &  +1   &  $-$0 &  $-$7 &  +6   &  $-$5 &  +4   & $-$11 & +10   &  $-$9 &  +8   & +15   & $-$14 & +13   & $-$12 \\
\hline
   4  &  +5   &  +6   &  +7   &  $-$0 &  $-$1 &  $-$2 &  $-$3 & $-$12 & $-$13 & $-$14 & $-$15 &  +8   &  +9   & +10   & +11   \\
   5  &  $-$4 &  +7   &  $-$6 &  +1   &  $-$0 &  +3   &  $-$2 & $-$13 & +12   & $-$15 & +14   &  $-$9 &  +8   & $-$11 & +10   \\
   6  &  $-$7 &  $-$4 &  +5   &  +2   &  $-$3 &  $-$0 &  +1   & $-$14 & +15   & +12   & $-$13 & $-$10 & +11   &  +8   &  $-$9 \\
   7  &  +6   &  $-$5 &  $-$4 &  +3   &  +2   &  $-$1 &  $-$0 & $-$15 & $-$14 & +13   & +12   & $-$11 & $-$10 &  +9   &  +8   \\
\hline   
	 8  &  +9   & +10   & +11   & +12   & +13   & +14   & +15   &  $-$0 &  $-$1 &  $-$2 &  $-$3 &  $-$4 &  $-$5 &  $-$6 &  $-$7 \\
   9  &  $-$8 & +11   & $-$10 & +13   & $-$12 & $-$15 & +14   &  +1   &  $-$0 &  +3   &  $-$2 &  +5   &  $-$4 &  $-$7 &  +6   \\
  10  & $-$11 &  $-$8 &  +9   & +14   & +15   & $-$12 & $-$13 &  +2   &  $-$3 &  $-$0 &  +1   &  +6   &  +7   &  $-$4 &  $-$5 \\
  11  & +10   &  $-$9 &  $-$8 & +15   & $-$14 & +13   & $-$12 &  +3   &  +2   &  $-$1 &  $-$0 &  +7   &  $-$6 &  +5   &  $-$4 \\
  12  & $-$13 & $-$14 & $-$15 &  $-$8 &  +9   & +10   & +11   &  +4   &  $-$5 &  $-$6 &  $-$7 &  $-$0 &  +1   &  +2   &  +3   \\
  13  & +12   & $-$15 & +14   &  $-$9 &  $-$8 & $-$11 & +10   &  +5   &  +4   &  $-$7 &  +6   &  $-$1 &  $-$0 &  $-$3 &  +2   \\
  14  & +15   & +12   & $-$13 & $-$10 & +11   &  $-$8 &  $-$9 &  +6   &  +7   &  +4   &  $-$5 &  $-$2 &  +3   &  $-$0 &  $-$1 \\
  15  & $-$14 & +13   & +12   & $-$11 & $-$10 &  +9   &  $-$8 &  +7   &  $-$6 &  +5   &  +4   &  $-$3 &  $-$2 &  +1   &  $-$0 \\
 \hline \hline
\end{tabular}%
}
\end{table}

\noindent
An inspection of this table yields as many as 35 distinguished triples, namely:
$$
\{1,2,3\}, \{1,4,5\}, \{1,6,7\}, \{1,8,9\}, \{1,10,11\}, \{1,12,13\}, \{1,14,15\},
$$
$$
\{2,4,6\}, \{2,5,7\}, \{2,8,10\}, \{2,9,11\}, \{2,12,14\}, \{2,13,15\},
 $$
$$
\{3,4,7\}, \{3,5,6\}, \{3,8,11\}, \{3,9,10\}, \{3,12,15\}, \{3,13,14\},
$$
$$
\{4,8,12\}, \{4,9,13\}, \{4,10,14\}, \{4,11,15\},
$$
$$
\{5,8,13\}, \{5,9,12\}, \{5,10,15\}, \{5,11,14\},
$$
$$
\{6,8,14\}, \{6,9,15\}, \{6,10,12\}, \{6,11,13\},
$$
$$
\{7,8,15\}, \{7,9,14\}, \{7,10,13\}, \{7,11,12\}.
$$

\begin{figure}[pth!]
	\centering
	\includegraphics[width=9truecm]{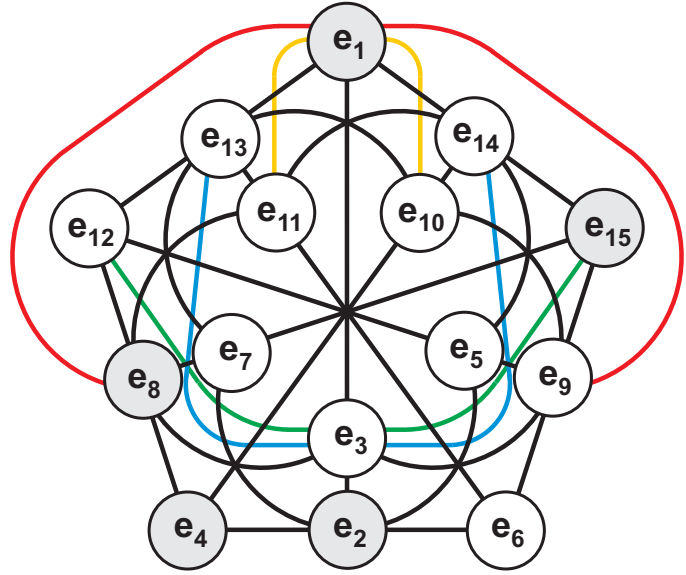}
	\caption{An illustration of the structure of PG$(3,2)$ that provides the multiplication law for sedenions. As in the previous case, the three imaginaries lying on the same line are such that the product of two of them yields the third one, sign disregarded.}
	\label{fgr1}
	\end{figure}

\noindent
Regarding the 15 imaginary units as points and the 35 distinguished trios of them as lines, we obtain a point-line incidence structure where each line has three points and each point is on seven lines, and which is isomorphic to PG$(3,2)$, the smallest projective space --- as depicted in Figure \ref{fgr1}. The latter figure employs a diagrammatical model of PG$(3,2)$ built, after Polster \cite{pol}, around the pentagonal model of the generalized quadrangle of type GQ$(2,2)$ whose 15 lines are illustrated by triples of points lying on black line-segments (10 of them) and/or black arcs of circles (5). The remaining 20 lines of PG$(3,2)$ comprise four distinct orbits: the yellow, red, blue and green one consisting, respectively, of the yellow ($\{1,10,11\}$), red
($\{1,8,9\}$), blue ($\{3,13,14\}$) and green ($\{3,12,15\}$) line and other four lines we get from each by rotation through 72 degrees around the center of the pentagon. 
It is not difficult to verify that we have 25 ordinary lines, namely
$$
\{1,2,3\}, \{1,4,5\}, \{1,6,7\}, \{1,8,9\}, \{1,10,11\}, \{1,12,13\}, \{1,14,15\},
$$
$$
\{2,4,6\}, \{2,5,7\}, \{2,8,10\}, \{2,9,11\}, \{2,12,14\}, \{2,13,15\},
 $$
$$
\{3,4,7\}, \{4,8,12\}, \{4,9,13\}, \{4,10,14\}, \{4,11,15\},
$$
$$
\{3,8,11\}, \{5,8,13\}, \{6,8,14\}, \{7,8,15\},
$$
$$
\{3,12,15\}, \{5,10,15\},  \{6,9,15\}, 
$$
and 10 defective ones, namely
$$
 \{3,5,6\},  \{3,9,10\}, \{3,13,14\},
$$
$$
\{5,9,12\}, \{5,11,14\},
$$
$$
\{6,10,12\}, \{6,11,13\},
$$
$$
\{7,9,14\}, \{7,10,13\}, \{7,11,12\}.
$$
Similarly, our sedenionic PG$(3,2)$ features two distinct types of points. 
A type-one point is such that four lines passing through it are ordinary, the remaining three being defective; such a point lies in the set 
$$\{3,5,6,7,9,10,11,12,13,14\} \equiv \alpha.$$
A type-two point is such that every line passing through it is ordinary; such a point belongs  to the set 
$$\{1,2,4,8,15\} \equiv \beta,$$
being illustrated by  gray shading in Figure \ref{fgr1}.
We see that all defective lines are of the same form, namely $\{\alpha, \alpha, \alpha \}$. The 25 ordinary lines split into two distinct families. Ten of them are of the form $\{\alpha, \beta, \beta\}$, namely
$$
\{1,2,3\}, \{1,4,5\}, \{1,8,9\}, \{1,14,15\}, 
$$
$$
\{2,4,6\}, \{2,8,10\}, \{2,13,15\},
 $$
$$
 \{4,8,12\},\{4,11,15\},
$$
$$
 \{7,8,15\},
$$
and the remaining 15 are of the form $\{\alpha, \alpha, \beta\}$, namely
$$
\{1,6,7\},  \{1,10,11\}, \{1,12,13\}, 
$$
$$
\{2,5,7\},  \{2,9,11\}, \{2,12,14\}, 
 $$
$$ 
\{3,4,7\}, \{4,9,13\}, \{4,10,14\}, 
$$
$$
\{3,8,11\}, \{5,8,13\}, \{6,8,14\}, 
$$
$$
\{3,12,15\}, \{5,10,15\},  \{6,9,15\}.
$$

\begin{figure}[pth!]
	\centering
	\includegraphics[width=4.5truecm]{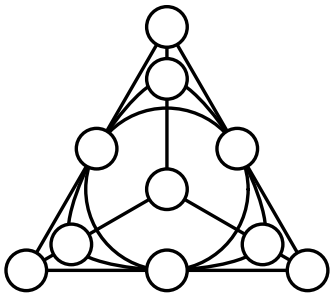}
	\caption{An illustrative portrayal of the Desargues configuration, built around the model of the Pasch configuration shown in Figure \ref{f2}: circles stand for  its points, whereas its lines are represented by triples of points on common straight segments (six), arcs of circles (three) and a big circle.}
\label{fgr2}
\end{figure}

\begin{figure}[pth!]
	\centering
	\psfrag{1}{\small 14}
	\psfrag{2}{\small 13}
	\psfrag{3}{\small 11}
	\psfrag{4}{\small 9}
	\psfrag{5}{\small 10}
	\psfrag{6}{\small 12}
	\psfrag{7}{\small 6}
	\psfrag{8}{\small 5}
	\psfrag{9}{\small 3}
	\psfrag{10}{\small ~7}
	\psfrag{11}{\small ~1}
	\psfrag{12}{\small ~2}
	\psfrag{13}{\small 4}
	\psfrag{14}{\small ~8}
	\psfrag{15}{\small 15}
	\includegraphics[width=12truecm]{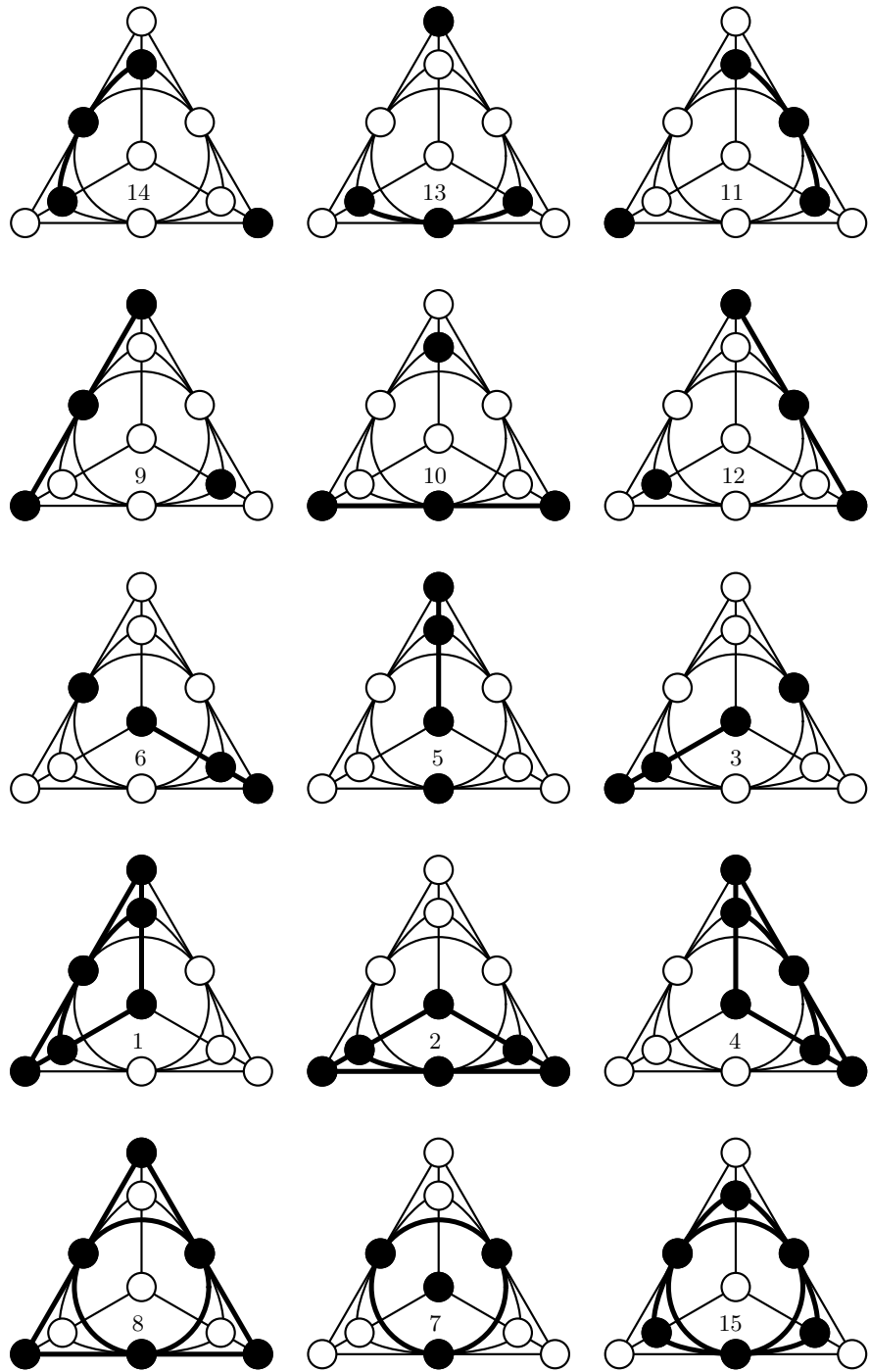}
	\caption{The fifteen geometric hyperplanes of the Desargues configuration. The hyperplanes are labelled by imaginary units of sedenions in such a way that --- as we shall verify in the next three figures --- the 35 lines of the Veldkamp space of the Desargues configuration are identical with the 35 distinguished triples of units, that is with the 35 lines of the PG$(3,2)$ shown in Figure \ref{fgr1}.}
\label{fgr3}	
\end{figure}

\begin{figure}[pth!]
	\centering
	\psfrag{1-4-10}{~7-9-14}
	\psfrag{2-5-10}{7-10-13}
	\psfrag{3-6-10}{7-11-12}
	\psfrag{1-2-9}{3-13-14}
	\psfrag{1-3-8}{5-11-14}
	\psfrag{2-3-7}{6-11-13}
	\psfrag{4-5-9}{3-9-10}
	\psfrag{4-6-8}{5-9-12}
	\psfrag{5-6-7}{6-10-12}
	\psfrag{7-8-9}{~3-5-6}
	\includegraphics[width=14truecm,clip=]{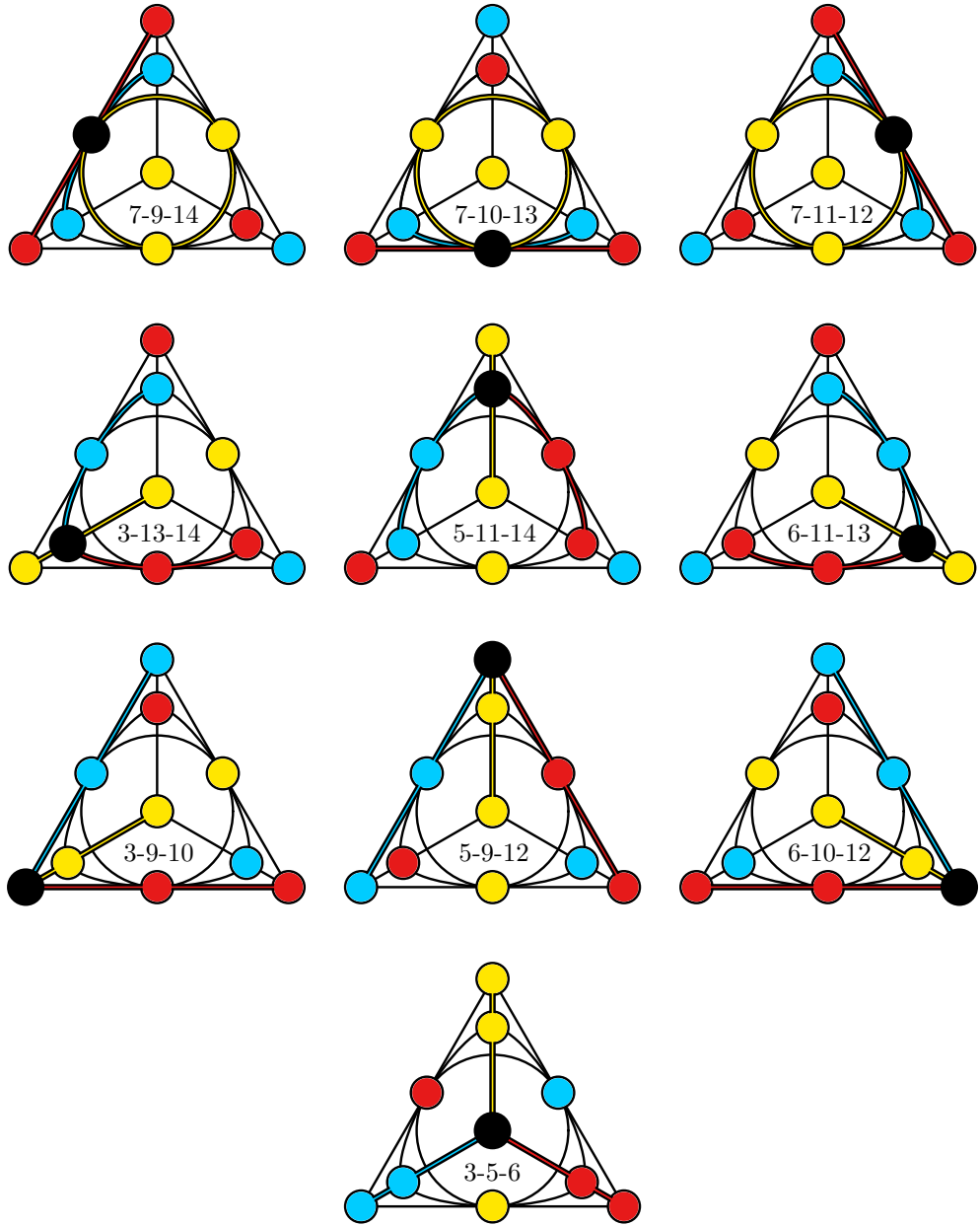}
	\caption{The ten Veldkamp lines of the Desargues configuration that represent the ten defective lines
	of the sedenionic PG$(3,2)$. Here, as well as in the next two figures, the three geometric hyperplanes 
	comprising a given Veldkamp line
	are distinguished by different colors, with their common elements (here just a single point) being colored black.
	For each Veldkamp line we also explicitly indicate its composition.}
	\label{fgr4}
\end{figure}

\begin{figure}[pth!]
	\centering
	\psfrag{7-12-13}{~~2-4-6}
	\psfrag{8-11-13}{~~1-4-5}
	\psfrag{9-11-12}{~~1-2-3}
	\psfrag{6-13-14}{~~4-8-12}
	\psfrag{5-12-14}{~~2-8-10}
	\psfrag{4-11-14}{~~1-8-9}
	\psfrag{3-13-15}{~~4-11-15}
	\psfrag{2-12-15}{~~2-13-15}
	\psfrag{1-11-15}{~~1-14-15}
	\psfrag{10-14-15}{~~7-8-15}
	\includegraphics[width=14truecm,clip=]{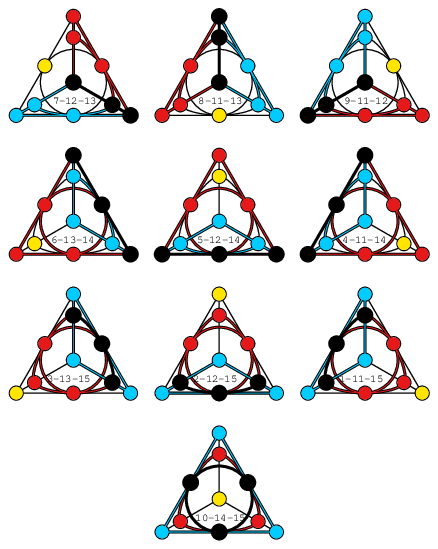}
	\caption{The ten Veldkamp lines of the Desargues configuration that represent the ten ordinary lines
	of the sedenionic PG$(3,2)$ of type $\{\alpha, \beta, \beta\}$.}
	\label{fgr5}
\end{figure}

\begin{figure}[pth!]
	\centering
	\psfrag{4-7-15}{~~6-9-15}
	\psfrag{5-8-15}{~~5-10-15}
	\psfrag{6-9-15}{~~3-12-15}
	\psfrag{1-7-14}{~~6-8-14}
	\psfrag{2-8-14}{~~5-8-13}
	\psfrag{3-9-14}{~~3-8-11}
	\psfrag{7-10-11}{~~1-6-7}
	\psfrag{8-10-12}{~~2-5-7}
	\psfrag{9-10-13}{~~3-4-7}
	\psfrag{1-6-12}{~2-12-14}
	\psfrag{2-4-13}{~~4-9-13}
	\psfrag{3-5-11}{~1-10-11}
	\psfrag{3-4-12}{~~2-9-11}
	\psfrag{2-6-11}{~1-12-13}
	\psfrag{1-5-13}{~4-10-14}
	\includegraphics[width=13truecm,clip=]{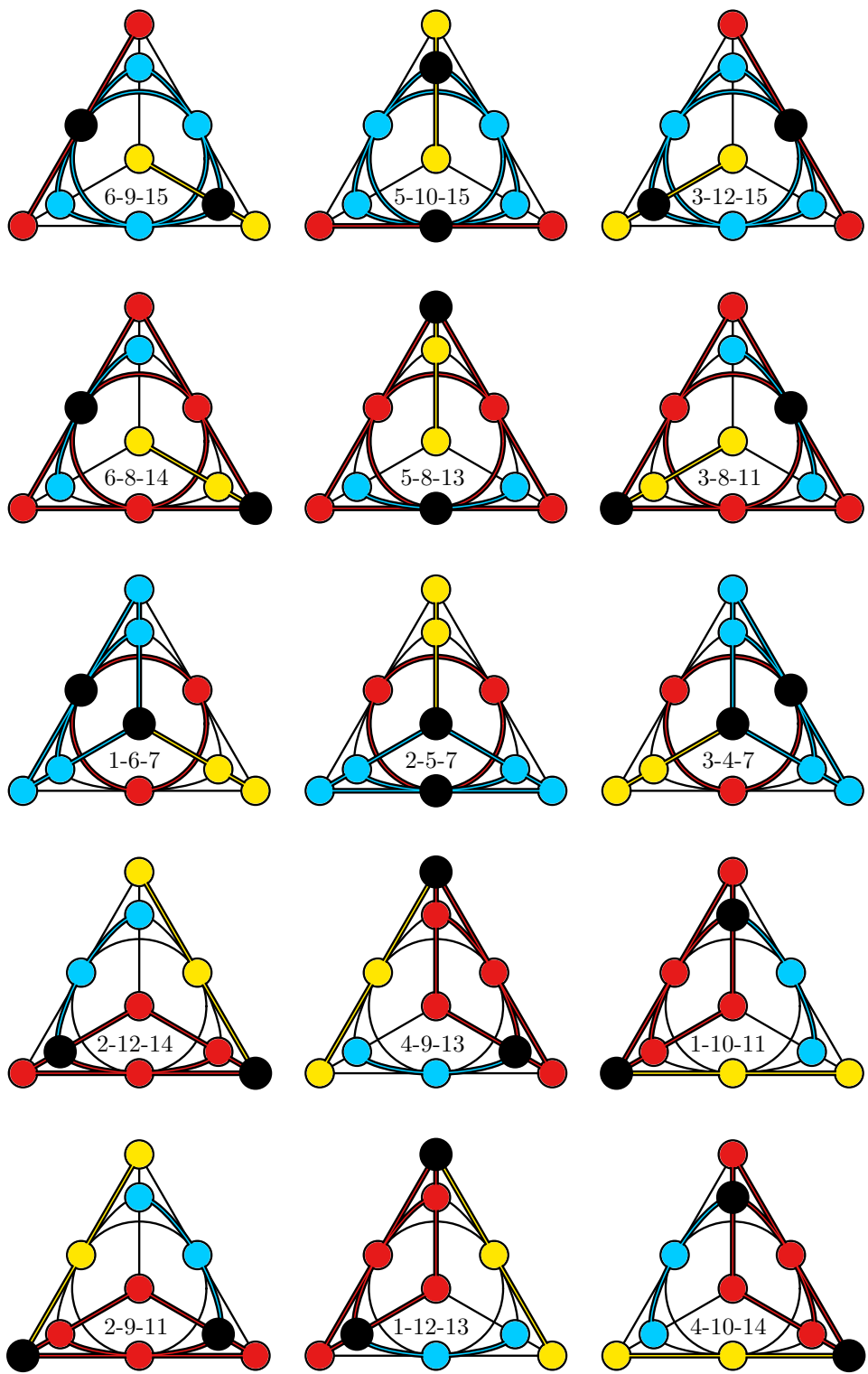}
	\caption{The fifteen Veldkamp lines of the Desargues configuration that represent the fifteen ordinary lines
	of the sedenionic PG$(3,2)$ of type $\{\alpha, \alpha, \beta\}$.}
	\label{fgr6}
\end{figure}

\begin{figure}[t]
	\centering
	\psfrag{1}{\tiny 1}
	\psfrag{2}{\tiny 2}
	\psfrag{3}{\tiny 3}
	\psfrag{4}{\tiny 4}
	\psfrag{5}{\tiny 5}
	\psfrag{6}{\tiny 6}
	\psfrag{7}{\tiny 7}
	\psfrag{8}{\tiny 8}
	\psfrag{9}{\tiny 9}
	\psfrag{10}{\tiny 10}
	\psfrag{11}{\tiny 11}
	\psfrag{12}{\tiny 12}
	\psfrag{13}{\tiny 13}
	\psfrag{14}{\tiny 14}
	\psfrag{15}{\tiny 15}
	\includegraphics[width=14truecm]{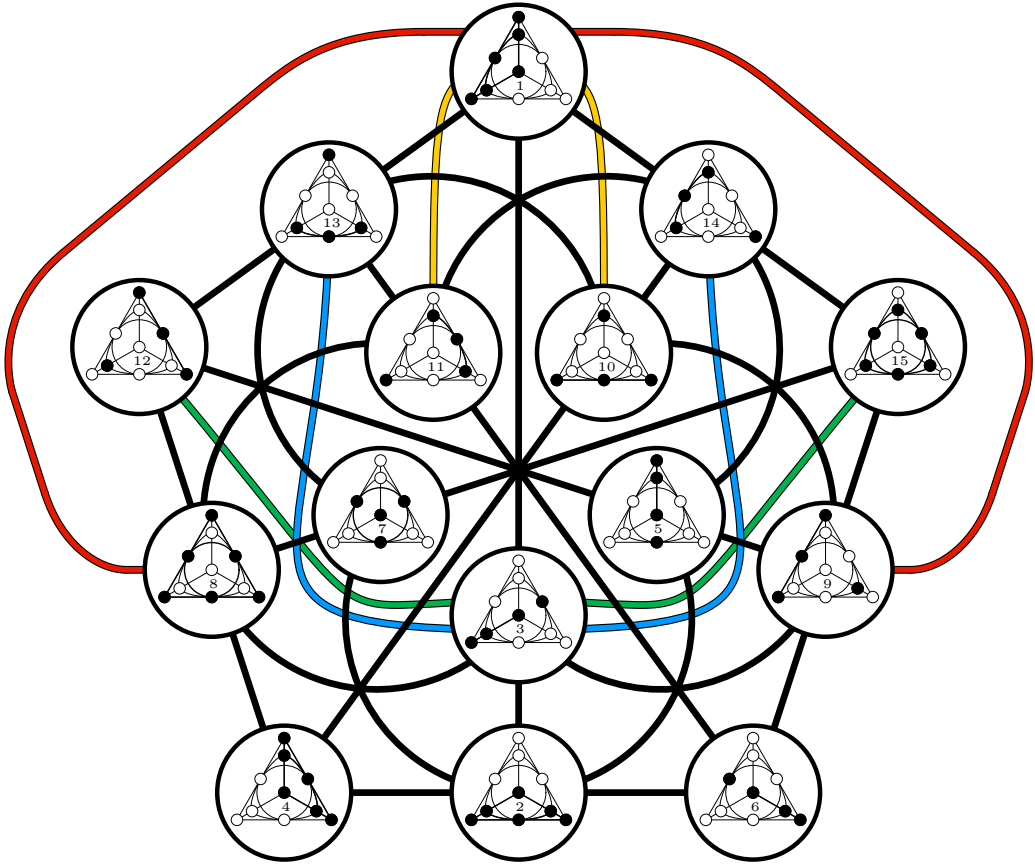}
	\caption{A compact graphical view of illustrating the bijection between 15 imaginary unit sedenions and 15 geometric hyperplanes of the Desargues configuration, as well as between 35 distinguished triples of units and 35 Veldkamp lines of the Desargues configuration. }
	\label{fgr7}
\end{figure}

A configuration $\mathcal{C}_4$ whose Veldkamp space reproduces the above-described partitions of points and lines of PG$(3,2)$ is, as demonstrated below,  the famous {\it Desargues} $(10_3)$-configuration, $\mathcal{D}$, which is one of the most prominent point-line incidence structures (see, e.\,g., \cite{piser}). Up to isomorphism, there exist altogether ten $(10_3)$-configurations. The Desargues configuration is, unlike the others, flag-transitive and the only one where for {\it each} of its points the three points that are not collinear with it lie on a line. This configuration, depicted in Figure \ref{fgr2} in a form showing its automorphism of order three, also lives in our sedenionic PG$(3,2)$; here, its points are the ten $\alpha$-points and its lines are all the defective lines.
In order to see that $\mathcal{V}(\mathcal{D}) \cong$ PG$(3,2)$ we shall first introduce, using our diagrammatical representation of $\mathcal{D}$, all 15 geometric hyperplanes of $\mathcal{D}$ --- Figure \ref{fgr3}. We see that they are indeed of two different forms, and  of required cardinality ten and five. A member of the former set comprises a point and three points not collinear with it. Such geometric hyperplane corresponds to a type-one (or $\alpha$-) point of PG$(3,2)$. A member of the latter set features six points located on four lines, with two lines per each point; this is nothing but the Pasch configuration we introduced in the previous section. Such a geometric hyperplane of $\mathcal{D}$ corresponds to a type-two (or $\beta$-) point of our PG$(3,2)$. It is also a straightforward task to verify that $\mathcal{V}(\mathcal{D})$ is endowed with 35 lines splitting into the required three families; those that correspond to defective lines of our sedenionic PG$(3,2)$ are shown in Figure \ref{fgr4}, while those that correspond to ordinary lines are depicted in Figure \ref{fgr5} (of type $\{\alpha, \beta, \beta\}$) and Figure \ref{fgr6} (of type $\{\alpha, \alpha, \beta\}$). Figure \ref{fgr7} offers a `condensed' view of the isomorphism $\mathcal{V}(\mathcal{D}) \cong$ PG$(3,2)$.

We shall finalize this section by pointing out that the existence of two different kinds of geometric hyperplanes of the Desargues configuration is closely connected with two well-known views of this configuration. The first one is as a pair of triangles that are in perspective from both a point and a line (Desargues' theorem), the point and the line forming a geometric hyperplane. The other view is  as the incidence sum of a complete quadrangle (i.\,e., a $(4_3,6_2)$-configuration) and a Pasch $(6_2,4_3)$-configuration \cite{bogepi}.

\newpage
\section{32-nions and the Cayley-Salmon $(15_4,20_3)$-configuration}
Our next case is $A_5$, or the 32-nions, whose multiplication properties are encoded in Table 3.
\begin{table}[pth!]
\centering
\caption{The multiplication table --- for better readability split into two parts --- of the imaginary unit 32-nions $e_a$, $1 \leq a \leq 31$, in short-handed notation $e_a \equiv a$ (and $e_0 \equiv 0$). } 
\vspace*{0.1cm}
\resizebox{\columnwidth}{!}{%
\begin{tabular}{||c|r|rr|rrrr|rrrrrrrrr||}
\hline \hline
   $*$   &  1    &  2    &  3    &  4    &  5    &  6    &  7    &  8    &  9    & 10    & 11    & 12    & 13    & 14    & 15  & ~~~~~~ \\
\hline
   1  &  $-$0 &  $-$3 &  +2   &  $-$5 &  +4   &  +7   &  $-$6 &  $-$9 &  +8   & +11   & $-$10 & +13   & $-$12 & $-$15 & +14  & \\
\hline
   2  &  +3   &  $-$0 &  $-$1 &  $-$6 &  $-$7 &  +4   &  +5   & $-$10 & $-$11 &  +8   &  +9   & +14   & +15   & $-$12 & $-$13 &\\
   3  &  $-$2 &  +1   &  $-$0 &  $-$7 &  +6   &  $-$5 &  +4   & $-$11 & +10   &  $-$9 &  +8   & +15   & $-$14 & +13   & $-$12 &\\
\hline
   4  &  +5   &  +6   &  +7   &  $-$0 &  $-$1 &  $-$2 &  $-$3 & $-$12 & $-$13 & $-$14 & $-$15 &  +8   &  +9   & +10   & +11   &\\
   5  &  $-$4 &  +7   &  $-$6 &  +1   &  $-$0 &  +3   &  $-$2 & $-$13 & +12   & $-$15 & +14   &  $-$9 &  +8   & $-$11 & +10   &\\
   6  &  $-$7 &  $-$4 &  +5   &  +2   &  $-$3 &  $-$0 &  +1   & $-$14 & +15   & +12   & $-$13 & $-$10 & +11   &  +8   &  $-$9 &\\
   7  &  +6   &  $-$5 &  $-$4 &  +3   &  +2   &  $-$1 &  $-$0 & $-$15 & $-$14 & +13   & +12   & $-$11 & $-$10 &  +9   &  +8   &\\
\hline   
	 8  &  +9   & +10   & +11   & +12   & +13   & +14   & +15   &  $-$0 &  $-$1 &  $-$2 &  $-$3 &  $-$4 &  $-$5 &  $-$6 &  $-$7 &\\
   9  &  $-$8 & +11   & $-$10 & +13   & $-$12 & $-$15 & +14   &  +1   &  $-$0 &  +3   &  $-$2 &  +5   &  $-$4 &  $-$7 &  +6   &\\
  10  & $-$11 &  $-$8 &  +9   & +14   & +15   & $-$12 & $-$13 &  +2   &  $-$3 &  $-$0 &  +1   &  +6   &  +7   &  $-$4 &  $-$5 &\\
  11  & +10   &  $-$9 &  $-$8 & +15   & $-$14 & +13   & $-$12 &  +3   &  +2   &  $-$1 &  $-$0 &  +7   &  $-$6 &  +5   &  $-$4 &\\
  12  & $-$13 & $-$14 & $-$15 &  $-$8 &  +9   & +10   & +11   &  +4   &  $-$5 &  $-$6 &  $-$7 &  $-$0 &  +1   &  +2   &  +3   &\\
  13  & +12   & $-$15 & +14   &  $-$9 &  $-$8 & $-$11 & +10   &  +5   &  +4   &  $-$7 &  +6   &  $-$1 &  $-$0 &  $-$3 &  +2   &\\
  14  & +15   & +12   & $-$13 & $-$10 & +11   &  $-$8 &  $-$9 &  +6   &  +7   &  +4   &  $-$5 &  $-$2 &  +3   &  $-$0 &  $-$1 &\\
  15  & $-$14 & +13   & +12   & $-$11 & $-$10 &  +9   &  $-$8 &  +7   &  $-$6 &  +5   &  +4   &  $-$3 &  $-$2 &  +1   &  $-$0 &\\
	\hline
  16  &  +17  &  +18  &  +19  &  +20  &  +21  &  +22  &  +23  &  +24  & +25   & +26   & +27   & +28   & +29   & +30   & +31   &\\
  17  & $-$16 &  +19  & $-$18 &  +21  & $-$20 & $-$23 &  +22  &  +25  & $-$24 & $-$27 & +26   & $-$29 & +28   & +31   & $-$30 &\\
  18  & $-$19 & $-$16 &  +17  &  +22  &  +23  & $-$20 & $-$21 &  +26  & +27   & $-$24 & $-$25 & $-$30 & $-$31 & +28   & +29   &\\
  19  &  +18  & $-$17 & $-$16 &  +23  & $-$22 &  +21  & $-$20 &  +27  & $-$26 & +25   & $-$24 & $-$31 & +30   & $-$29 & +28   &\\
  20  & $-$21 & $-$22 & $-$23 & $-$16 &  +17  &  +18  &  +19  &  +28  & +29   & +30   & +31   & $-$24 & $-$25 & $-$26 & $-$27  & \\
  21  &  +20  & $-$23 & +22   & $-$17 & $-$16 & $-$19 &  +18  &  +29  & $-$28 & +31   & $-$30 &  +25  & $-$24 & +27   &  $-$26 &\\
  22  &  +23  & +20   & $-$21 & $-$18 &  +19  & $-$16 & $-$17 &  +30  & $-$31 & $-$28 & +29   &  +26  & $-$27 & $-$24 &  +25   &\\ 
	23  & $-$22 & +21   & +20   & $-$19 & $-$18 & +17   & $-$16 &  +31  &  +30  & $-$29 & $-$28 &  +27  &  +26  & $-$25 &  $-$24 &\\
  24  & $-$25 & $-$26 & $-$27 & $-$28 & $-$29 & $-$30 & $-$31 & $-$16 &  +17  & +18   & +19   &  +20  &  +21  &  +22  &  +23   &\\
  25  & +24   & $-$27 & +26   & $-$29 & +28   & +31   & $-$30 & $-$17 & $-$16 & $-$19 & +18   & $-$21 &  +20  &  +23  &  $-$22 &\\
  26  & +27   &  +24  & $-$25 & $-$30 & $-$31 & +28   & +29   & $-$18 &  +19  & $-$16 & $-$17 & $-$22 & $-$23 &  +20  &  +21 &\\
  27  & $-$26 &  +25  & +24   & $-$31 &  +30  & $-$29 & +28   & $-$19 & $-$18 &  +17  & $-$16 & $-$23 &  +22  & $-$21 &  +20  & \\
  28  & +29   &  +30  & +31   &  +24  & $-$25 & $-$26 & $-$27 & $-$20 &  +21  &  +22  &  +23  & $-$16 & $-$17 & $-$18 &  $-$19 &  \\
  29  & $-$28 & +31   & $-$30 &  +25  & +24   &  +27  & $-$26 & $-$21 & $-$20 &  +23  & $-$22 &  +17  & $-$16 &  +19  &  $-$18 &\\
  30  & $-$31 & $-$28 & +29   &  +26  & $-$27 &  +24  &  +25  & $-$22 & $-$23 & $-$20 &  +21  &  +18  & $-$19 & $-$16 &  +17 &\\
	31  & +30   & $-$29 & $-$28 &  +27  & +26   & $-$25 &  +24  & $-$23 &  +22  & $-$21 & $-$20 &  +19  &  +18  & $-$17 &  $-$16 &\\
 \hline 
\end{tabular}%
}
\begin{tabular}{cc}
& \\
\end{tabular}
\resizebox{\columnwidth}{!}{%
\begin{tabular}{||c|rrrrrrrrrrrrrrrr||}
\hline 
$*$  &  16   &  17   &  18   &  19   &  20   &  21   &  22   &  23   &  24   &  25   &  26   &  27   &   28  &  29   &  30   & 31 \\
\hline
  1  & $-$17 &  +16  &  +19  & $-$18 &  +21  & $-$20 & $-$23 &  +22  & +25   & $-$24 & $-$27 & +26   & $-$29 &  +28  &  +31  & $-$30 \\
	\hline
  2  & $-$18 & $-$19 &  +16  &  +17  &  +22  &  +23  & $-$20 & $-$21 & +26   &  +27  & $-$24 & $-$25 & $-$30 & $-$31 &  +28  & +29 \\
  3  & $-$19 &  +18  & $-$17 &  +16  &  +23  & $-$22 &  +21  & $-$20 & +27   & $-$26 &  +25  & $-$24 & $-$31 & +30   & $-$29 & +28 \\
	\hline
  4  & $-$20 & $-$21 & $-$22 & $-$23 &  +16  &  +17  &  +18  &  +19  & +28   & +29   & +30   &  +31  & $-$24 & $-$25 & $-$26 & $-$27 \\
  5  & $-$21 &  +20  & $-$23 &  +22  & $-$17 &  +16  & $-$19 &  +18  & +29   & $-$28 & +31   & $-$30 &  +25  & $-$24 & +27   & $-$26 \\
  6  & $-$22 &  +23  &  +20  & $-$21 & $-$18 &  +19  &  +16  & $-$17 & +30   & $-$31 & $-$28 &  +29  & +26   & $-$27 & $-$24 & +25 \\
  7  & $-$23 & $-$22 &  +21  &  +20  & $-$19 & $-$18 &  +17  &  +16  & +31   &  +30  & $-$29 & $-$28 &  +27  &  +26  & $-$25 & $-$24 \\ 
\hline
	8  & $-$24 & $-$25 & $-$26 & $-$27 & $-$28 & $-$29 & $-$30 & $-$31 & +16   &  +17  &  +18  &  +19  &  +20  &  +21  &  +22  & +23 \\
  9  & $-$25 & +24   & $-$27 & +26   & $-$29 &  +28  & +31   & $-$30 & $-$17 &  +16  & $-$19 &  +18  & $-$21 &  +20  &  +23  & $-$ 22 \\
 10  & $-$26 & +27   &  +24  & $-$25 & $-$30 & $-$31 & +28   &  +29  & $-$18 &  +19  &  +16  & $-$17 & $-$22 & $-$23 &  +20  & +21 \\
 11  & $-$27 & $-$26 &  +25  & +24   & $-$31 & +30   & $-$29 &  +28  & $-$19 & $-$18 &  +17  &  +16  & $-$23 &  +22  & $-$21 & +20
\\
 12  & $-$28 & +29   & +30   & +31   &  +24  & $-$25 & $-$26 & $-$27 & $-$20 &  +21  &  +22  &  +23  &  +16  & $-$17 & $-$18 & $-$19 \\
 13  & $-$29 & $-$28 & +31   & $-$30 &  +25  &  +24  & +27   & $-$26 & $-$21 & $-$20 &  +23  & $-$22 &  +17  & +16   &  +19  & $-$18 \\
 14  & $-$30 & $-$31 & $-$28 & +29   &  +26  & $-$27 & +24   &  +25  & $-$22 & $-$23 & $-$20 &  +21  &  +18  & $-$19 &  +16  & +17 \\
 15  & $-$31 & +30   & $-$29 & $-$28 &  +27  &  +26  & $-$25 &  +24  & $-$23 &  +22  & $-$21 & $-$20 &  +19 &  +18   & $-$17 & +16
\\
\hline
 16  & $-$0  & $-$1  & $-$2  & $-$3  & $-$4  & $-$5  & $-$6  & $-$7  & $-$8  & $-$9  & $-$10 & $-$11 & $-$12 & $-$13 & $-$14 & $-$15 \\
 17  &  +1   & $-$0  &  +3   & $-$2  &  +5   & $-$4  & $-$7  &  +6   & +9    & $-$8  & $-$11 &  +10  & $-$13 & +12   & +15   & $-$14 \\
 18  &  +2   & $-$3  & $-$0  &  +1   &  +6   &  +7   & $-$4  & $-$5  & +10   & +11   & $-$8  & $-$9  & $-$14 & $-$15 & +12   & +13 \\
 19  &  +3   &  +2   & $-$1  & $-$0  &  +7   & $-$6  &  +5   & $-$4  & +11   & $-$10 &  +9   & $-$8  & $-$15 & +14   & $-$13 & +12 \\
 20  &  +4   & $-$5  & $-$6  & $-$7  & $-$0  &  +1   &  +2   &  +3   & +12   & +13   & +14   &  +15  & $-$8  & $-$9  & $-$10 & $-$11 \\
 21  &  +5   & +4    & $-$7  &  +6   & $-$1  & $-$0  & $-$3  &  +2   & +13   & $-$12 & +15   & $-$14 &  +9   & $-$8  &  +11  & $-$10 \\
 22  &  +6   & +7    & +4    & $-$5  & $-$2  & +3    & $-$0  & $-$1  & +14   & $-$15 & $-$12 &  +13  &  +10  & $-$11 & $-$8  & +9 \\ 
 23  &  +7   & $-$6  & +5    &  +4   & $-$3  & $-$2  &  +1   & $-$0  &  +15  & +14   & $-$13 & $-$12 &  +11  &  +10  & $-$9  & $-$8\\
 24  &  +8   & $-$9  & $-$10 & $-$11 & $-$12 & $-$13 & $-$14 & $-$15 & $-$0  &  +1   & +2    &  +3   &  +4   &  +5   &  +6   & +7 \\
 25  &  +9   &  +8   & $-$11 &  +10  & $-$13 & +12   & +15   & $-$14 & $-$1  & $-$0  & $-$3  &  +2   & $-$5  &  +4   &  +7   & $-$6\\
 26  &  +10  &  +11  & +8    & $-$9  & $-$14 & $-$15 & +12   &  +13  & $-$2  &  +3   & $-$0  & $-$1  & $-$6  & $-$7  &  +4   & +5 \\
 27  &  +11  & $-$10 & +9    &  +8   & $-$15 &  +14  & $-$13 &  +12  & $-$3  & $-$2  &  +1   & $-$0  & $-$7  &  +6   & $-$5  & +4 \\
 28  &  +12  &  +13  & +14   &  +15  & +8    & $-$9  & $-$10 & $-$11 & $-$4  &  +5   &  +6   &  +7   & $-$0  & $-$1  & $-$2  & $-$3 \\
 29  &  +13  & $-$12 & +15   & $-$14 & +9    &  +8   &  +11  & $-$10 & $-$5  & $-$4  &  +7   & $-$6  &  +1   & $-$0  &  +3   & $-$2 \\
 30  &  +14  & $-$15 & $-$12 &  +13  & +10   & $-$11 &  +8   & +9    & $-$6  & $-$7  & $-$4  &  +5   &  +2   & $-$3  & $-$0  & +1 \\
 31  &  +15  &  +14  & $-$13 & $-$12 & +11   &  +10  & $-$9  & +8    & $-$7  &  +6   & $-$5  & $-$4  &  +3   &  +2   & $-$1  & $-$0 \\
 \hline \hline
\end{tabular}%
}
\end{table}
From this table we infer the existence of 155 distinguished triples of imaginary units. Regarding the 31 imaginary units of $A_5$ as points and the 155 distinguished triples of them as lines, we obtain a point-line incidence structure where each line has three points and each point is on 15 lines, and which is isomorphic to PG$(4,2)$. We next find that 65 lines of this space are defective and 90 ordinary. However, unlike the preceding two cases, there are {\it three} different types of points in our 32-nionic PG$(4,2)$:
ten $\alpha$-points, $$\alpha \equiv \{7,11,13,14,19,21,22,25,26,28\},$$ each of which is on nine defective and six ordinary lines; fifteen $\beta$-points, $$\beta \equiv \{3,5,6,9,10,12,15,17,18,20,23,24,27,29,30\},$$ each of which is on seven defective and eight ordinary lines; and six $\gamma$-points, $$\gamma \equiv \{1,2,4,8,16,31\},$$ each of them being on fifteen ordinary (and, hence, on zero defective) lines. This stratification of the point-set of PG$(4,2)$ leads, in turn, to two different kinds of defective lines and three distinct kinds of ordinary lines. 
Concerning the former, there are 45 of them of type $\{\alpha,\alpha,\beta\}$, namely
$$
\{3,13,14\}, \{3,21,22\}, \{3,25,26\}, 
$$
$$
\{5,11,14\}, \{5,19,22\}, \{5,25,28\},
$$ 
$$
\{6,11,13\}, \{6,19,21\}, \{6,26,28\}, 
$$
$$
\{7,9,14\}, \{7,10,13\}, \{7,11,12\}, \{7,17,22\}, \{7,18,21\}, \{7,19,20\}, 
$$ $$
\{7,25,30\}, \{7,26,29\}, \{7,27,28\}, 
$$
$$
\{9,19,26\}, \{9,21,28\}, 
$$
$$
\{10,19,25\}, \{10,22,28\}, 
$$
$$
\{11,17,26\}, \{11,18,25\}, \{11,19,24\}, \{11,21,30\}, \{11,22,29\}, \{11,23,28\}, 
$$
$$
\{12,21,25\}, \{12,22,26\}, 
$$
$$
\{13,17,28\}, \{13,19,30\}, \{13,20,25\}, \{13,21,24\}, \{13,22,27\}, \{13,23,26\},
$$
$$
\{14,18,28\}, \{14,19,29\}, \{14,20,26\}, \{14,21,27\}, \{14,22,24\}, \{14,23,25\}, 
$$
$$
\{15,19,28\}, \{15,21,26\}, \{15,22,25\}, 
$$
and 20 of type $\{\beta,\beta,\beta\}$, namely
$$
\{3,5,6\}, \{3,9,10\}, \{3,17,18\}, \{3,29,30\}, 
$$ 
$$
\{5,9,12\}, \{5,17,20\}, \{5,27,30\}, 
$$
$$
\{6,10,12\}, \{6,18,20\}, \{6,27,29\},
$$
$$ 
\{9,17,24\}, \{9,23,30\}, 
$$
$$
\{10,18,24\}, \{10,23,29\}, 
$$ 
$$
\{12,20,24\}, \{12,23,27\}, 
$$ 
$$
\{15,17,30\}, \{15,18,29\}, \{15,20,27\}, \{15,23,24\}.
$$ 
As per the latter, one finds 15 of them of type  $\{\beta,\beta,\beta\}$, namely
$$
\{3,12,15\}, \{3,20,23\}, \{3,24,27\}, 
$$
$$
\{5,10,15\}, \{5,18,23\}, \{5,24,29\}, 
$$
$$
\{6,9,15\}, \{6,17,23\}, \{6,24,30\}, 
$$
$$
\{9,18,27\}, \{9,20,29\}, 
$$
$$
\{10,17,27\}, \{10,20,30\}, 
$$
$$
\{12,17,29\}, \{12,18,30\}, 
$$
60 of type  $\{\alpha,\beta,\gamma \}$, namely
$$
\{1, 6, 7\}, \{1, 10, 11\}, \{1, 12, 13\}, \{1, 14, 15\}, \{1, 18, 19\}, \{1, 20, 21\}, 
$$ $$
\{1, 22, 23\}, \{1, 24, 25\}, 
\{1, 26, 27\}, \{1, 28, 29\}, 
$$ 
$$
\{2, 5, 7\}, \{2, 9, 11\}, \{2, 12, 14\}, \{2, 13, 15\}, \{2, 17, 19\}, \{2, 20, 22\}, 
$$ $$
\{2, 21, 23\}, \{2, 24, 26\}, 
\{2, 25, 27\}, \{2, 28, 30\}, 
$$
$$
\{3, 4, 7\}, \{3, 8, 11\}, \{3, 16, 19\}, \{3, 28, 31\}, 
$$ $$
\{4, 9, 13\}, \{4, 10, 14\}, \{4, 11, 15\}, \{4, 17, 21\}, \{4, 18, 22\}, 
$$ $$
\{4, 19, 23\}, 
\{4, 24, 28\}, \{4, 25, 29\}, 
\{4, 26, 30\}, 
$$
$$
\{5, 8, 13\}, \{5, 16, 21\}, \{5, 26, 31\}, 
$$ $$
\{6, 8, 14\}, \{6, 16, 22\}, \{6, 25, 31\}, 
$$ $$
\{7, 8, 15\}, \{7, 16, 23\}, \{7, 24, 31\}, 
$$ $$
\{8, 17, 25\}, \{8, 18, 26\}, \{8, 19, 27\}, \{8, 20, 28\}, \{8, 21, 29\}, \{8, 22, 30\}, 
$$ $$
\{9, 16, 25\}, \{9, 22, 31\}, 
$$ $$
\{10, 16, 26\}, \{10, 21, 31\}, 
$$ $$
\{11, 16, 27\}, \{11, 20, 31\}, 
$$ $$
\{12, 16, 28\}, \{12, 19, 31\}, 
$$ $$
\{13, 16, 29\}, \{13, 18, 31\}, 
$$ $$
\{14, 16, 30\}, \{14, 17, 31\}, 
$$
and, finally, 15 of type $\{\beta, \gamma, \gamma \}$, namely

$$
\{1, 2, 3\}, \{1, 4, 5\}, \{1, 8, 9\}, \{1, 16, 17\}, \{1, 30, 31\}, 
$$ $$
\{2, 4, 6\}, \{2, 8, 10\}, \{2, 16, 18\}, \{2, 29, 31\}, 
$$ $$
\{4, 8, 12\}, \{4, 16, 20\}, \{4, 27, 31\}, 
$$ $$
\{8, 16, 24\}, \{8, 23, 31\}, 
$$ $$
\{5, 16, 31\}.
$$

\begin{figure}[pth!]
	\centering
	\includegraphics[width=9truecm]{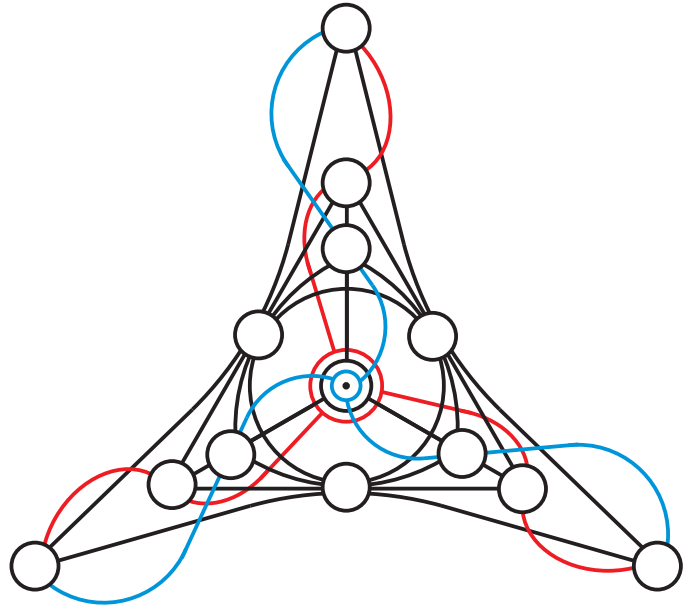}
	\caption{An illustration of the structure of the $(15_4,20_3)$-configuration, built around the model of the Desargues configuration shown in Figure \ref{fgr2}. The five points added to the Desargues configuration are the three peripheral points and the red and blue point in the center. The ten lines added are three lines denoted by red color, three blue lines, three lines joining pairwise the three peripheral points and the line that comprises the three points in the center of the figure, that is the ones represented by a bigger red circle, a smaller blue circle and a medium-sized black one.}
	\label{fig1}
	\end{figure}

\begin{figure}[pth!]
	\centering
	\includegraphics[width=14truecm]{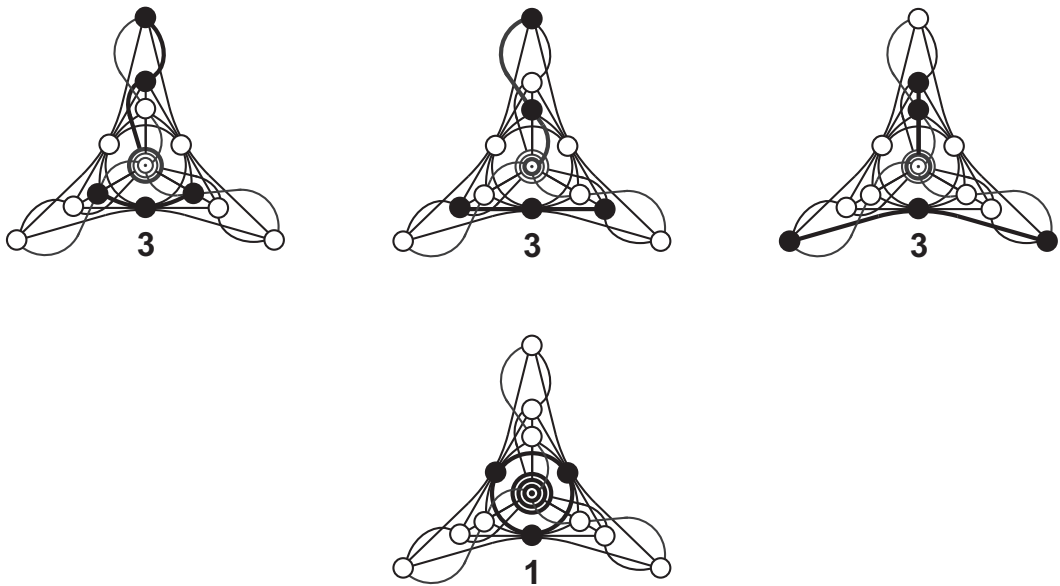}
	\caption{The ten geometric hyperplanes of the $(15_4,20_3)$-configuration of type one; the number below a subfigure indicates how many hyperplane's copies we get by rotating
	the particular subfigure through 120 degrees around its center.}
	\label{fig2}
	\end{figure}

\begin{figure}[pth!]
	\centering
	\includegraphics[width=14truecm]{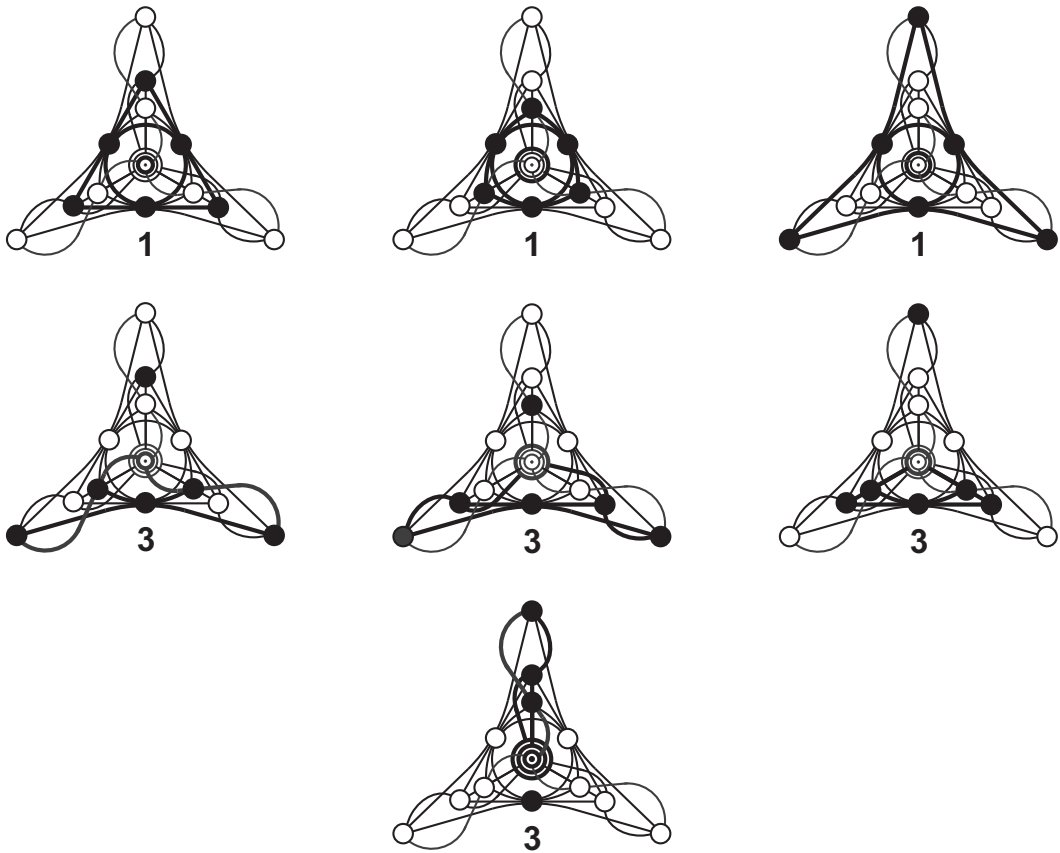}
	\caption{The fifteen geometric hyperplanes of the $(15_4,20_3)$-configuration of type two.}
	\label{fig3}
	\end{figure}

\begin{figure}[pth!]
	\centering
	\includegraphics[width=14truecm]{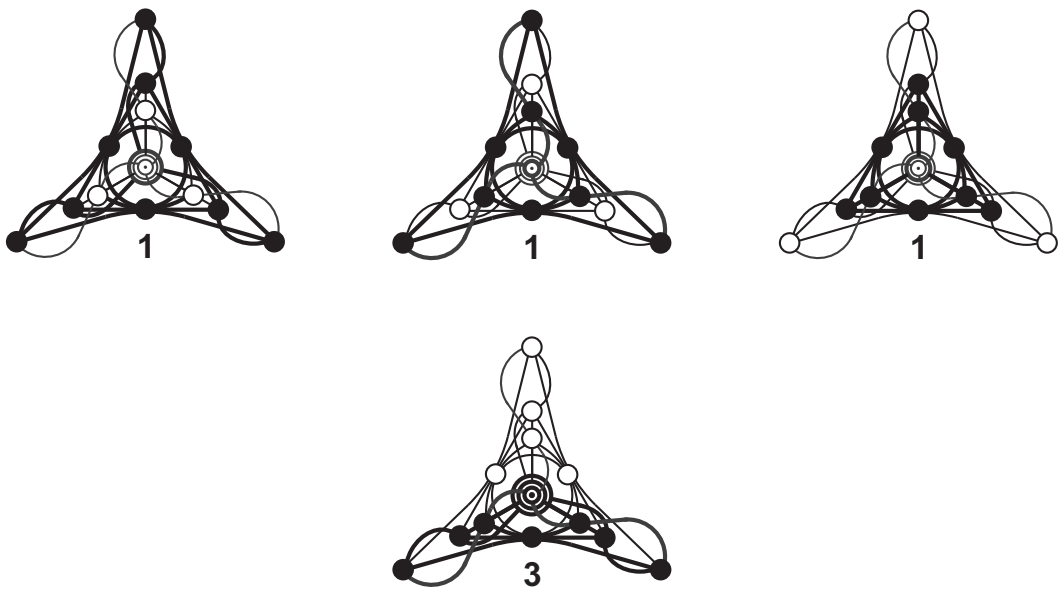}
	\caption{The six geometric hyperplanes of the $(15_4,20_3)$-configuration of type three.}
	\label{fig4}
	\end{figure}

A point-line configuration $\mathcal{C}_5$ whose Veldkamp space accounts for these stratifications of both the point- and line-set of our 32-nionic PG$(4,2)$ is of type $(15_4,20_3)$.\footnote{It is worth mentioning here that there exists another remarkable configuration whose Veldkamp space does the same job for us, namely the generalized quadrangle of order two, also known as the Cremona-Richmond $(15_3)$-configuration (see, e.\,g., \cite{twoq}).
However, this configuration is {\it triangle-free} and so it can$not$ contain the Desargues configuration as dictated by the nesting property of the Cayley-Dickson algebras.} This configuration is formed within our PG$(4,2)$ by 15 $\beta$-points and 20 defective lines of $\{\beta,\beta,\beta\}$ type and its structure is sketched in Figure \ref{fig1}. It is a rather easy task to verify that this particular $(15_4,20_3)$-configuration possesses 31 distinct geometric hyperplanes that fall into three different types. A type-one hyperplane consists of a pair of skew lines at maximum distance from each other;
there are, as depicted in Figure \ref{fig2},  ten hyperplanes of this type and they correspond to $\alpha$-points of PG$(4,2)$. A type-two hyperplane features a point and all the points not collinear with it, the latter forming --- not surprisingly --- the Pasch configuration; there are, as shown in 
Figure \ref{fig3}, fifteen hyperplanes of this type and their counterparts are $\beta$-points of PG$(4,2)$. A type-three hyperplane is identical with the Desargues configuration; we find, as portrayed in Figure \ref{fig4}, altogether six guys of this type, each standing for a $\gamma$-point of PG$(4,2)$.

		\begin{figure}[pth!]
	\centering
	\includegraphics[width=14truecm]{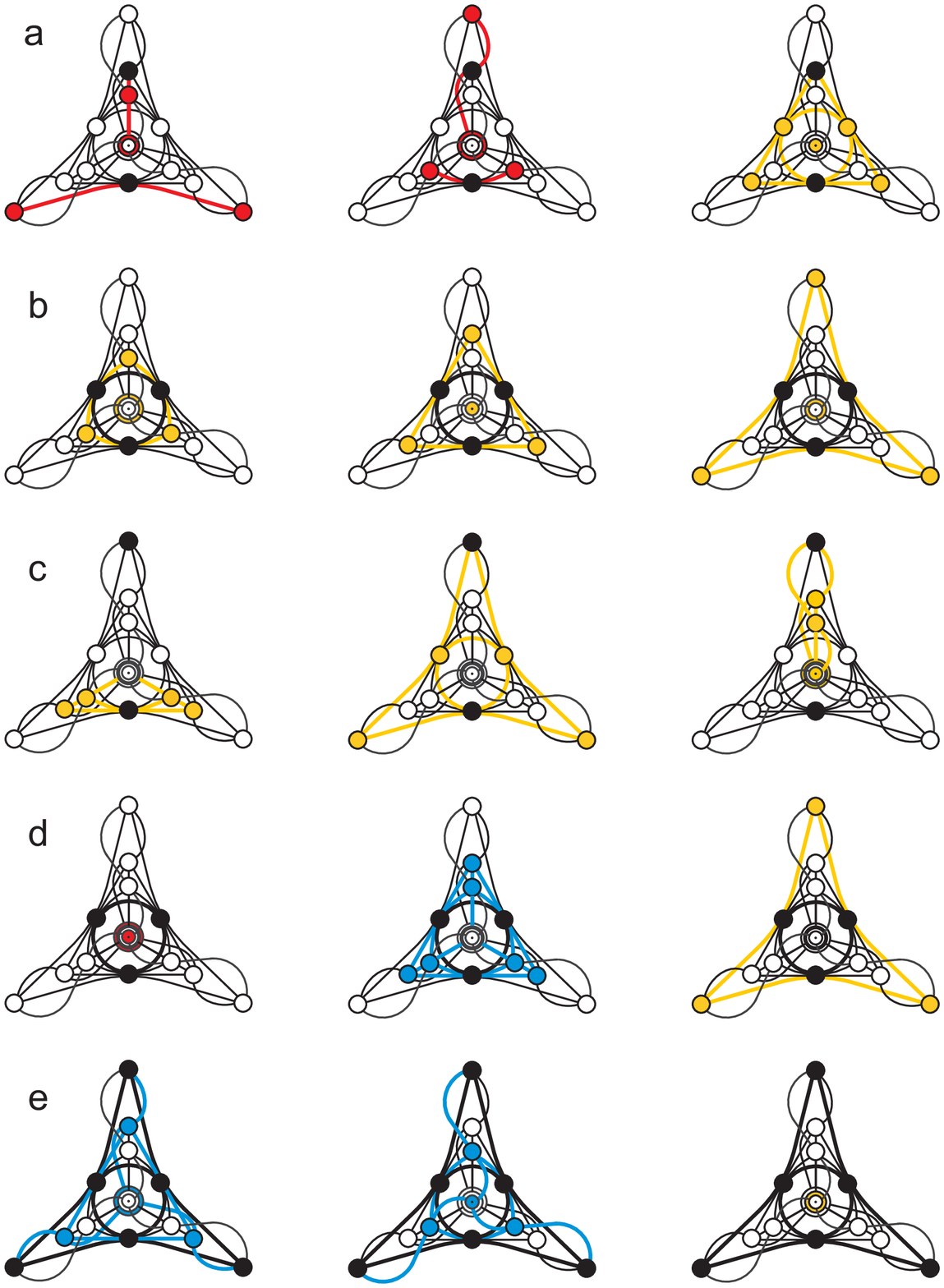}
	\caption{The five types of Veldkamp lines of the $(15_4,20_3)$-configuration. Here, unlike Figures 8 to 10, each representative  of a geometric hyperplane is drawn separately and different colors are used to distinguish between different hyperplane types: red is reserved for type one, yellow for type two and blue for type three hyperplanes. As before, black color denotes the core of a Veldkamp line, that is the elements common to all the three hyperplanes comprising it.}
	\label{fig5}
	\end{figure}
	
We also find that our $(15_4,20_3)$-configuration yields 155 Veldkamp lines that are, as expected,
of five different types. A type-I Veldkamp line, shown in Figure \ref{fig5}a, features two hyperplanes of type one and a type-two hyperplane and its core consists of two points that are at maximum distance from each other; there are ${10 \choose 2} = 15 \times 6/2 = 45$ Veldkamp lines of this type and they correspond to defective lines of PG$(4,2)$ of type $\{\alpha,\alpha,\beta\}$. A type-II Veldkamp line, featured in Figure \ref{fig5}b, is composed of three hyperplanes of type two that share three points on a common line; there are, obviously, 20 Veldkamp lines of this type, having for their counterparts defective lines of PG$(4,2)$ of type $\{\beta, \beta, \beta\}$. A type-III Veldkamp line, portrayed in Figure \ref{fig5}c, also consists of three hyperplanes of type two, but in this case the three common points are pairwise at maximum distance from each other; a quick count leads to 15 Veldkamp lines of this type, these being in a bijection with 15 ordinary lines of PG$(4,2)$ of type $\{\beta, \beta, \beta\}$. Next, it comes a type-IV Veldkamp line, depicted in Figure \ref{fig5}d, which exhibits a hyperplane of each type and whose core is composed of a line and a point at the maximum distance from it;
since for each line of our $(15_4,20_3)$-configuration there are three points at maximum distance from it, there are $20 \times 3 = 60$ Veldkamp lines of this type, having their twins in ordinary lines of PG$(4,2)$ of type $\{\alpha,\beta,\gamma \}$. Finally, we meet a type-V Veldkamp line, sketched in Figure \ref{fig5}e, which is endowed with two hyperplanes of type three and a single one of type two, and whose core is isomorphic to the Pasch configuration; hence, we have ${6 \choose 2} = 15$ Veldkamp lines of this type, being all representatives of ordinary lines of PG$(4,2)$ of type $\{\beta, \gamma, \gamma \}$.

Before embarking on the final case to be dealt with in detail, it is worth having a closer look at our $(15_4,20_3)$-configuration and pointing out its intimate relation with famous Pascal's Mystic Hexagram. If six arbitrary points are chosen on a conic section and joined by line segments in any order to form a hexagon, then the three pairs of opposite sides of the hexagon meet in three points that lie on a straight line, the latter being called the Pascal line. Taking the permutations of the six points, one obtains 60 different hexagons. Thus, the so-called complete Pascal hexagon determines altogether 60 Pascal lines, which generate a remarkable configuration of 146 points and 110 lines called the {\it hexagrammum mysticum}, or the complete Pascal figure
(for the most comprehensive, applet-based representation of this remarkable geometrical object, see \cite{norma}).
Both the points and lines of the complete Pascal figure split into several distinct families, usually named after their discoverers in the first half of the 19-th century. We are concerned here with the 15 Salmon points  and the 20 Cayley lines (see, e.\,g. \cite{lord,cory}) which form a $(15_4,20_3)$-configuration.  This configuration is discussed in some detail in \cite{bogepi}, where it is also depicted (Figure 6) and called the {\it Cayley-Salmon} $(15_4,20_3)$-configuration. And it is precisely this Cayley-Salmon $(15_4,20_3)$-configuration which our 32-nionic $(15_4,20_3)$-configuration is isomorphic to. The same configuration is also portrayed in Figure 8 of \cite{norma}. In the latter work, two different views/interpretations of the configuration are also mentioned. One is as three pairwise-disjoint triangles that are in perspective from a line, in which case the centers of perspectivity are guaranteed by Desargues' theorem to also lie on a line; we just stress here that these two lines form a geometric hyperplane (of type one, see Figure \ref{fig2}). The other view of the figure takes any point of the configuration to be the center of perspectivity of two quadrangles whose six pairs of corresponding sides meet necessarily in the points of a Pasch configuration; again, the point and the associated Pasch configuration form a geometric hyperplane (of type two, see Figure \ref{fig3}). Obviously, we can offer one more view of the configuration, that stemming from the existence of type-three hyperplanes, namely as the incidence sum of a Desargues configuration and three triangles on a commmon side (see Figure \ref{fig4}).

\section{64-nions and a $(21_5,35_3)$-configuration}
The final algebra we shall treat in sufficient detail is $A_6$, or the 64-nions.
From the corresponding multiplication table, which due to its size we do not show here but which is freely available at {\tt http://jjj.de/tmp-zero-divisors/mult-table-64-ions.txt}, we infer the existence of 651 distinguished triples of imaginary units. Regarding the 63 imaginary units of 64-nions as points and the 651 distinguished triples of them as lines, we obtain a point-line incidence structure where each line has three points and each point is on 31 lines, and which is isomorphic to PG$(5,2)$. Following the usual procedure, we find that 350 lines of this space are defective and 301 ordinary. Likewise the preceding case, we encounter {\it three} different types of points in our 64-nionic PG$(5,2)$:
35 $\alpha$-points, each of which is on 21 defective and 10 ordinary lines; 21 $\beta$-points, each of which is on 15 defective and 16 ordinary lines; and seven $\gamma$-points,  each of them being on 31 ordinary (and, hence, on zero defective) lines. This stratification of the point-set of PG$(5,2)$ leads, in turn, to three different kinds of defective lines and four distinct kinds of ordinary lines. Out of 350 defective lines, we find 105 of type $\{\alpha, \alpha, \alpha\}$, 210 of type $\{\alpha, \alpha, \beta\}$ and 35 of type $\{\beta, \beta, \beta\}$. On the other hand, 301 ordinary lines are partitioned into 105 guys of type $\{\alpha, \beta, \beta\}$, 70 of type $\{\alpha, \alpha, \gamma\}$, 105 of type $\{\alpha, \beta, \gamma\}$ and 21 of type $\{\beta, \gamma, \gamma\}$. 

\begin{figure}[t]
	
	\centerline{\includegraphics[width=4truecm]{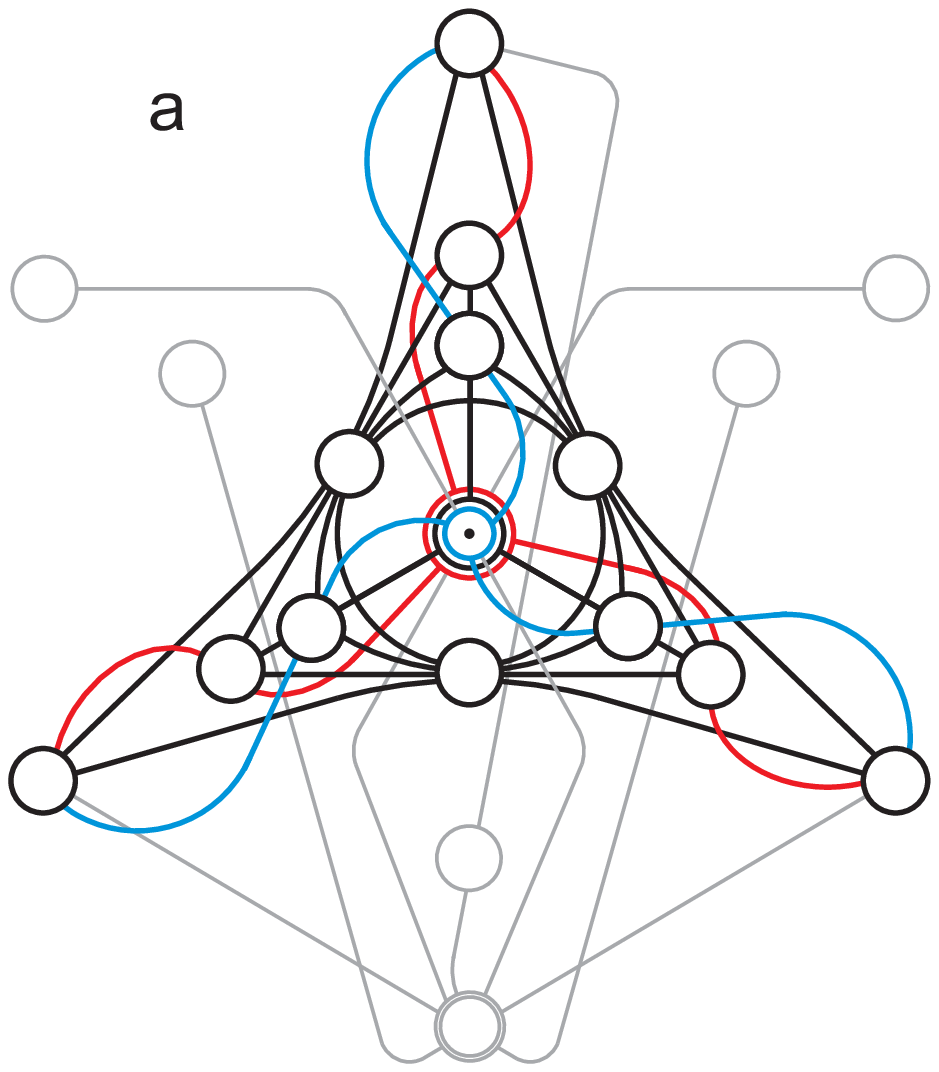}\hspace*{.6cm}\includegraphics[width=4truecm]{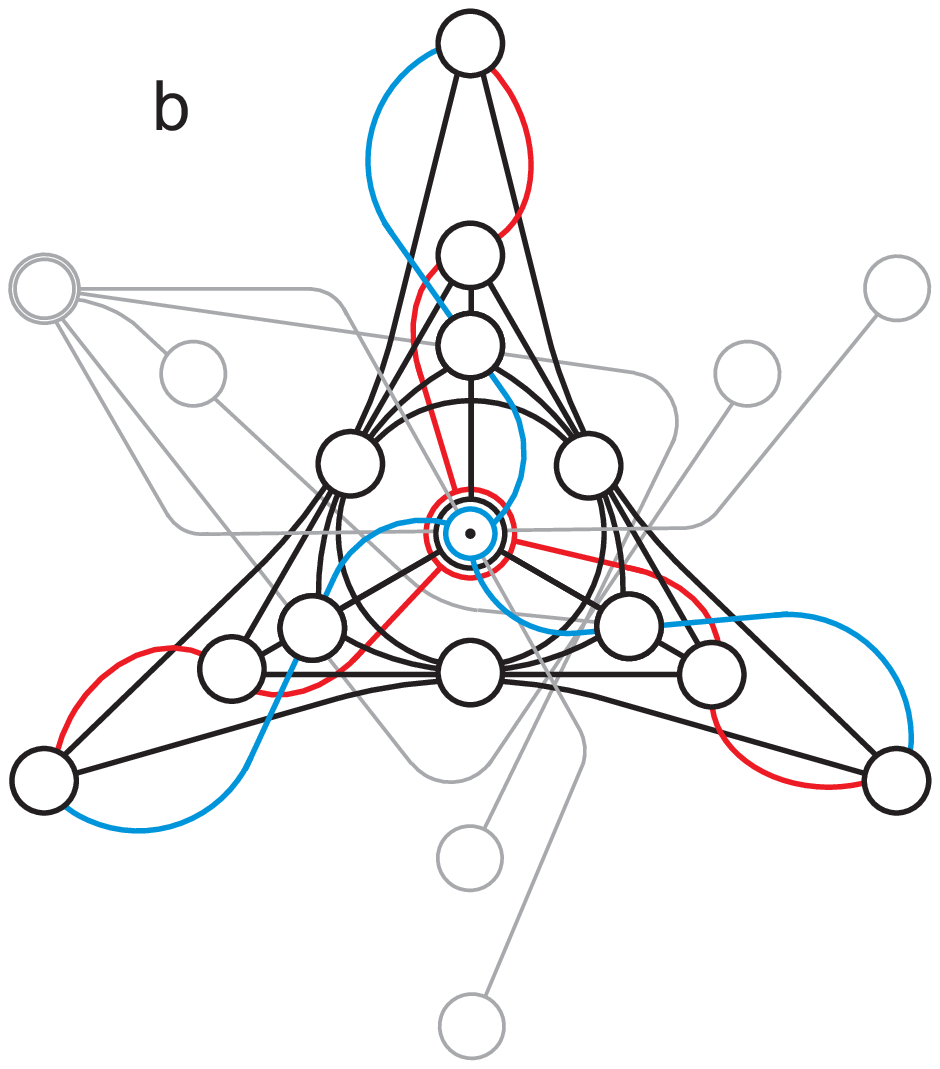}\hspace*{.6cm}\includegraphics[width=4truecm]{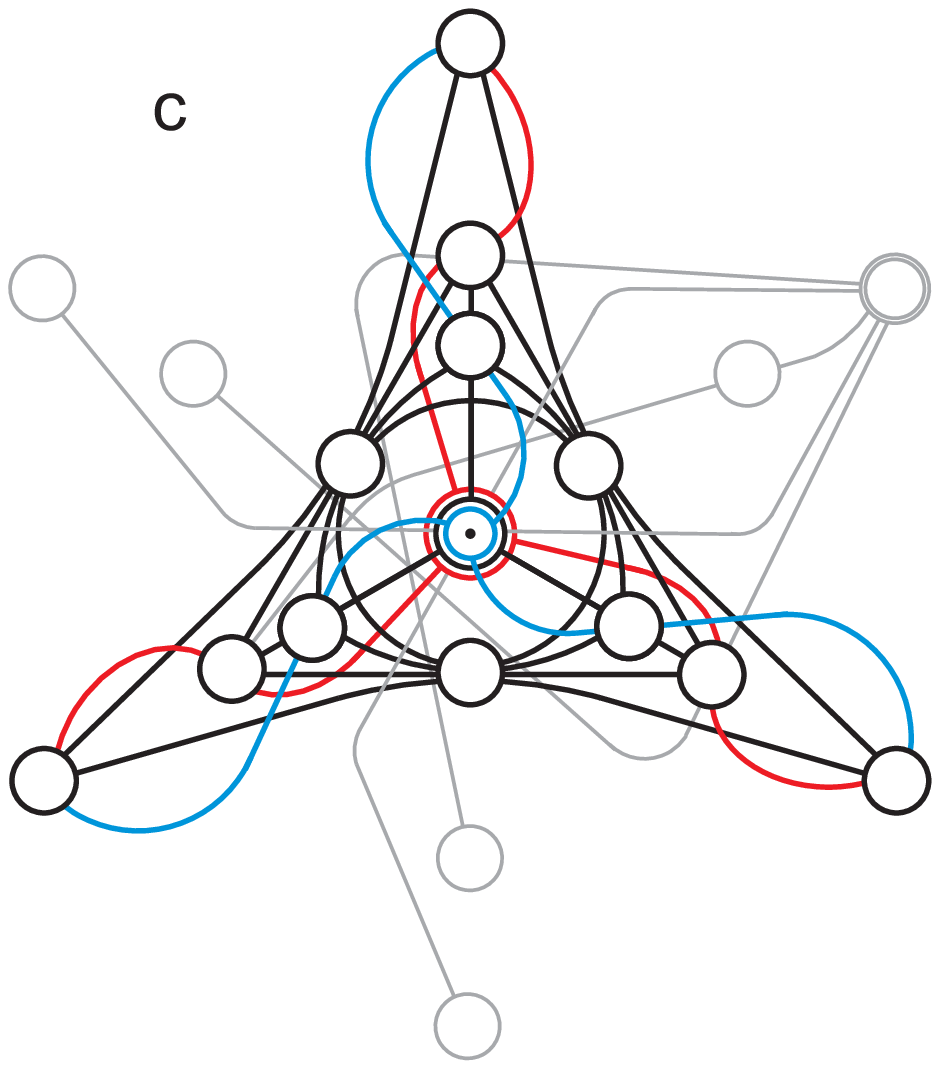}}
	\vspace*{.5cm}
	\centerline{\includegraphics[width=4.4truecm]{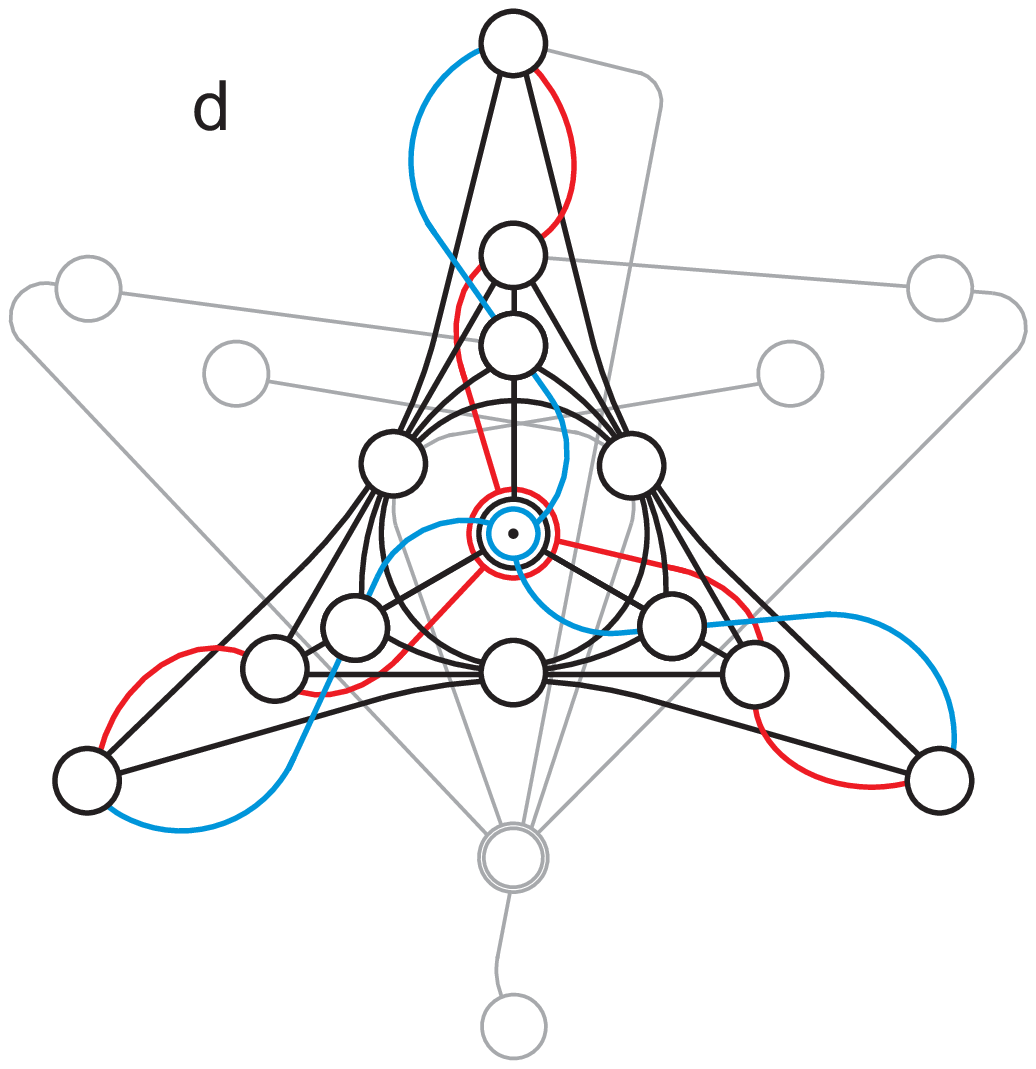}\hspace*{.6cm}\includegraphics[width=4truecm]{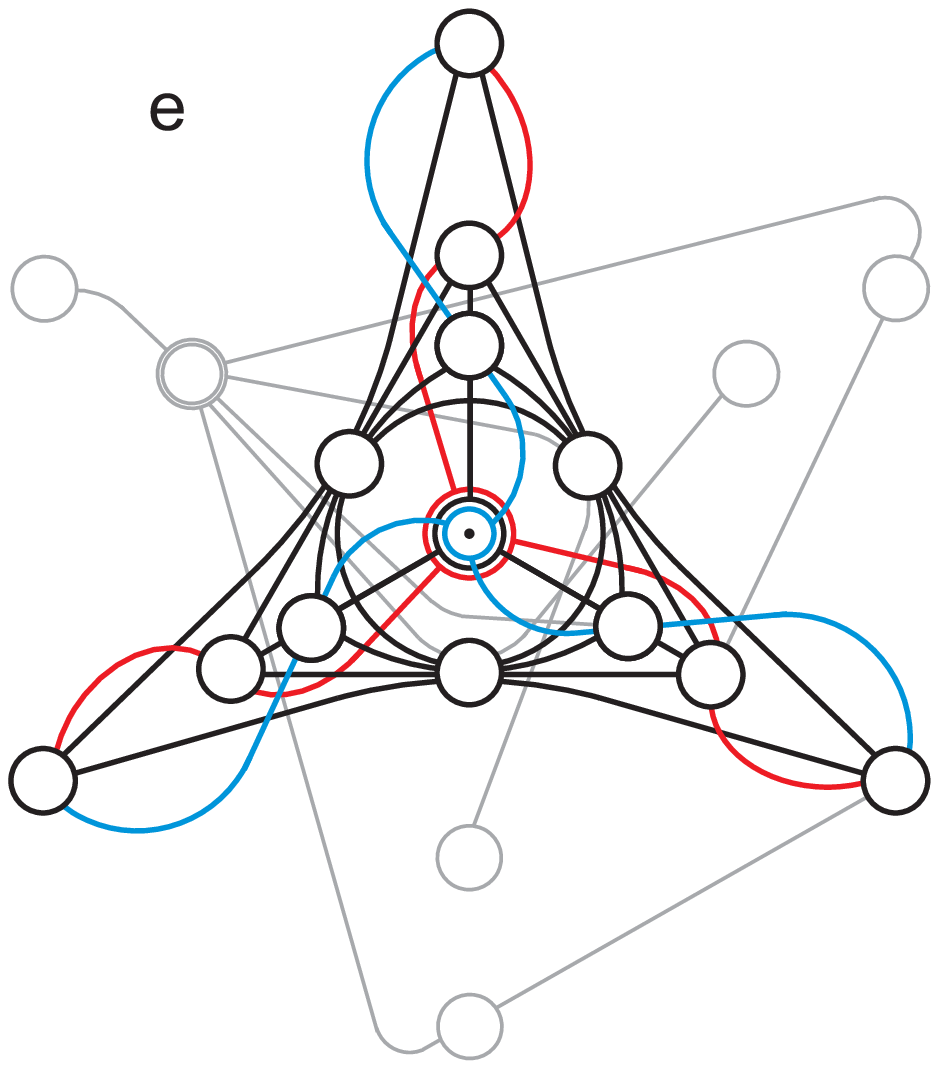}\hspace*{.6cm}\includegraphics[width=4truecm]{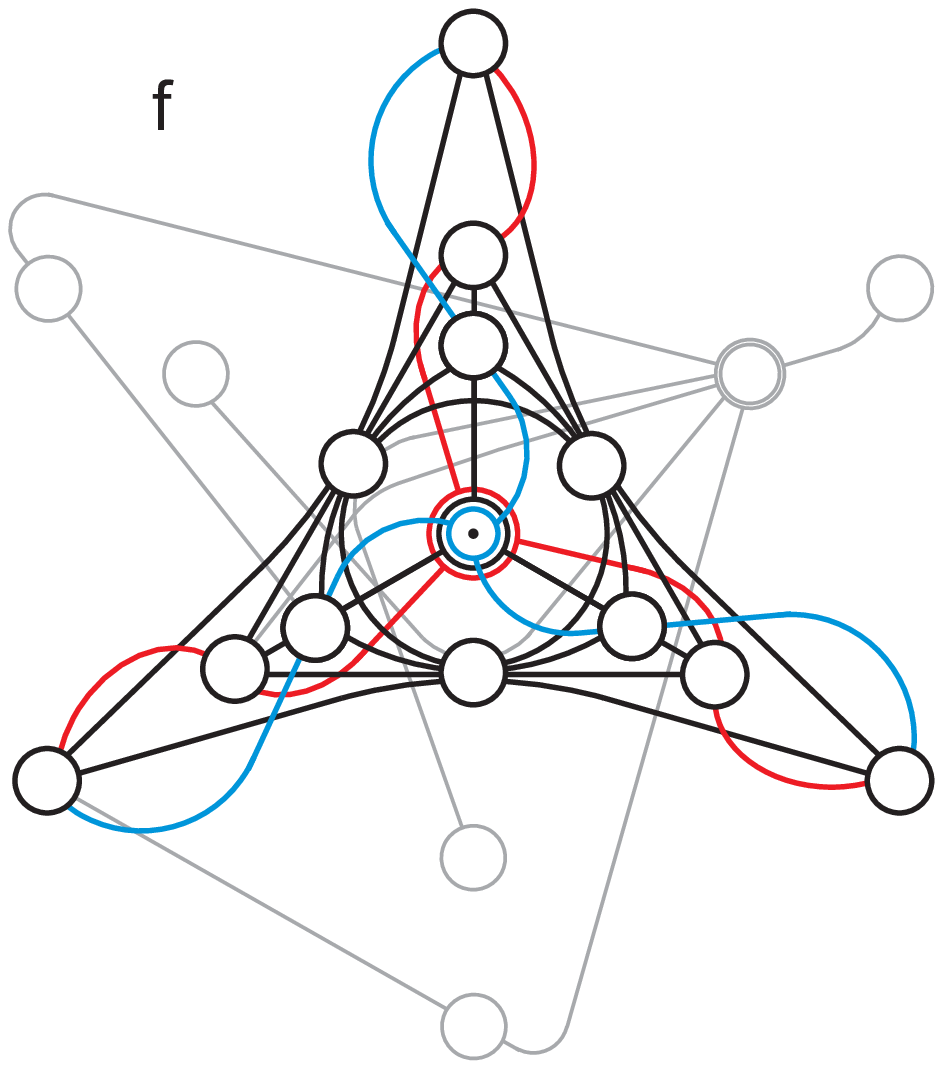}}
	\caption{An illustration of the structure of the $(21_5,35_3)$-configuration, built around the model of the Cayley-Salmon $(15_4,20_3)$-configuration shown in Figure \ref{fig1}.}
	\label{f64}
	\end{figure}

\pagebreak
The Veldkamp space mimicking such a fine structure of PG$(5,2)$ is that of a particular
$(21_5,35_3)$-configuration, $\mathcal{C}_6$, that also lives in our PG$(5,2)$ and whose points are the 21 $\beta$-points and whose lines are the 35 defective lines of $\{\beta, \beta, \beta\}$ type. To visualise this configuration, we build it around the model of the Cayley-Salmon $(15_4,20_3)$-configuration of 32-nions shown in Figure \ref{fig1}. Given the Cayley-Salmon configuration, there are six points and 15 lines to be added to yield our $(21_5,35_3)$-configuration, and this is to be done in such a way that the configuration we started with forms a geometric hyperplane in it. As putting all the lines into a single figure would make the latter look rather messy, in Figure \ref{f64} we briefly illustrate this construction by drawing six different figures, each featuring all six additional points (gray) but only five out of 15 additional lines (these lines being also drawn in gray color), namely those passing through a selected additional point (represented by a doubled circle).
Employing this handy diagrammatical representation, one can verify that our $(21_5,35_3)$-configuration exhibits 63 geometric hyperplanes that fall into three distinct types. A type-one hyperplane consists of a line and its complement, which is the Pasch configuration; there are 35 distinct hyperplanes of this form, each corresponding to an $\alpha$-point of our PG$(5,2)$. A type-two hyperplane comprises a point and its complement, which is the Desargues configuration; there are 21 hyperplanes of this form, each having a $\beta$-point for its PG$(5,2)$ counterpart. Finally,
a type-three hyperplane is isomorphic to the Cayley-Salmon configuration; there are seven distinct guys of this type, each answering to a 	
$\gamma$-point of the PG$(5,2)$. We leave it with the interested reader to verify by themselves that the Veldkamp space of our $(21_5,35_3)$-configuration indeed features 651 lines that do fall into the above-mentioned seven distinct kinds.

\begin{figure}[t]
	\centering
	\includegraphics[width=6truecm, angle=90]{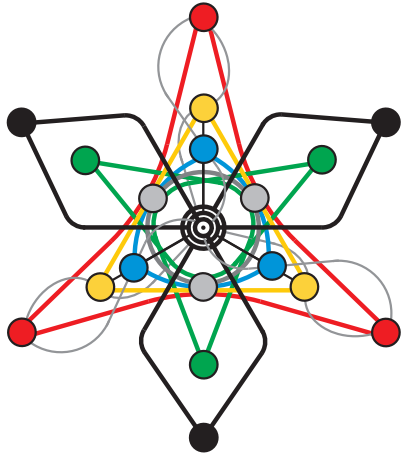}
	\caption{A `generalized Desargues' view of the $(21_5,35_3)$-configuration.}
	\label{f64view}
	\end{figure}

As in the previous two cases, we shall briefly describe a couple of interesting views of our $(21_5,35_3)$-configuration, both related to type-one hyperplanes.
The first one is as four triangles in perspective from a line where the points of perspectivity of six pairs of
them form a Pasch configuration, the line and the Pasch configuration comprising a geometric hyperplane (compare with the first view of both the Desargues and the Cayley-Salmon configuration). This is sketched in Figure
\ref{f64view}, where the four triangles are denoted, in boldfacing, by green, red, yellow and blue color, the
line of perspectivity by boldfaced gray color, and the points of perspectivity of
pairs of triangles (together with the corresponding lines they lie on and that are also boldfaced) by black color.
The other view is as three complete quadrangles that are pairwise in perspective in such a way that the three points of perspectivity lie on a line and where the six triples of their corresponding sides meet at points located on a Pasch configuration, again the line and the Pasch configuration forming a geometric hyperplane (compare with the second view of both the Desargues and the Cayley-Salmon configuration).

\section{$2^N$-nions and a $\left({N+1 \choose 2}_{N-1}, {N+1 \choose 3}_{3}\right)$-configuration}
\begin{figure}[t]
	\centering
	\includegraphics[width=14truecm]{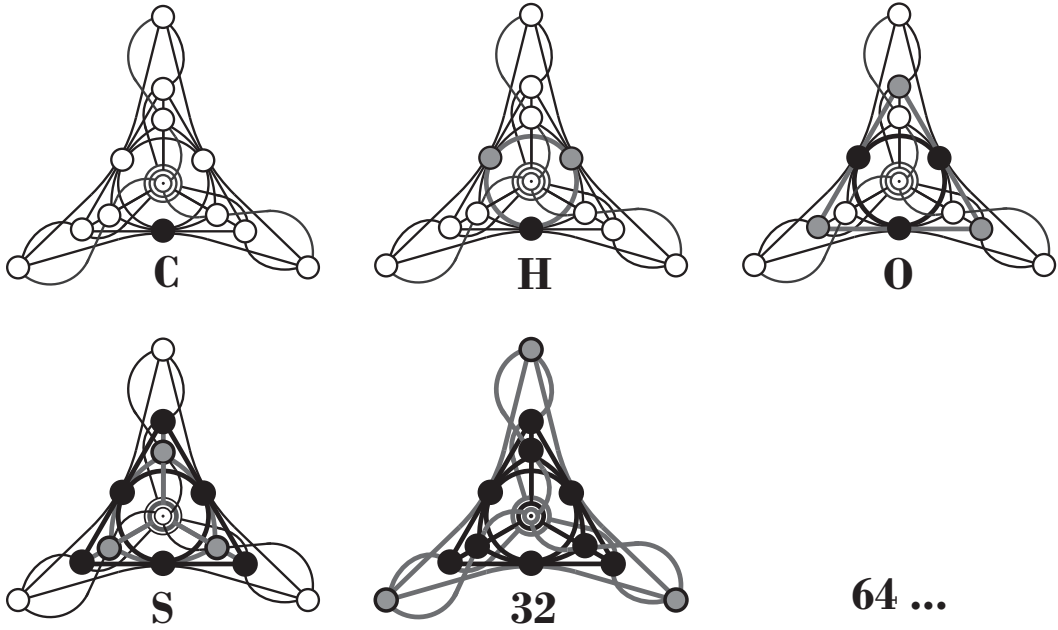}
	\caption{A nested hierarchy of finite $\left({N+1 \choose 2}_{N-1}, {N+1 \choose 3}_{3}\right)$-configurations of $2^N$-nions for $1 \leq N \leq 5$ when embedded in the Cayley-Salmon configuration ($N=5$).}
	\label{nested2}
	\end{figure}

At this point it is quite easy to spot the general pattern. If one also includes
the trivial cases of complex numbers ($N=1$), where the relevant geometry is just a single point i.\,e. the $(1_0,0_3)$-configuration, and quaternions ($N=2$), whose geometry is a single line i.\,e. the $(3_1,1_3)$-configuration, we obtain the following nested sequence of configurations
whose Veldkamp spaces capture the stratification/partition of the point- and line-sets of the $2^N$-nionic PG$(N-1,2)$, $N$ being a positive integer, 
$$(1_0,0_3),$$ $$(3_1,1_3),$$ $$(6_2,4_3),$$ $$(10_3,10_3),$$ $$(15_4,20_3),$$ $$(21_5,35_3),$$
$$\ldots,$$ 
$$\left({N+1 \choose 2}_{N-1}, {N+1 \choose 3}_{3}\right),$$
$$\ldots.$$
It is curious to notice that the first entry represents a {\it triangular} number, while the second one is a {\it tetrahedral} number, or triangular pyramidal number. In other words, we get a nested sequence of {\it binomial}  $({r+k-1 \choose r}_r, {r+k-1 \choose k}_k)$-configurations with $r = N-1$ and $k=3$, whose properties have very recently been discussed in a couple of interesting papers \cite{prap,pepr}.  The first few configurations are shown, in a form where the configurations are {\it nested} inside each other, in Figure \ref{nested2}.

A particular character of this nesting is reflected in the structure of geometric hyperplanes. Denoting our generic  $\left({N+1 \choose 2}_{N-1}, {N+1 \choose 3}_{3}\right)$-configuration by ${\cal C}_N$, we can express the types of geometric hyperplanes of the above-discussed cases in a compact form as follows

\bigskip
\begin{tabular}{llllll}
${\cal C}_1$: & $\varnothing$,&                                 &                                 &  & \\
${\cal C}_2$: & ${\cal C}_1$, &                                 &                                 &  & \\ 
${\cal C}_3$: & ${\cal C}_2$, & ${\cal C}_1 \sqcup {\cal C}_1$, &                                 &  & \\
${\cal C}_4$: & ${\cal C}_3$, & ${\cal C}_2 \sqcup {\cal C}_1$, &                                 &  & \\
${\cal C}_5$: & ${\cal C}_4$, & ${\cal C}_3 \sqcup {\cal C}_1$, & ${\cal C}_2 \sqcup {\cal C}_2$, &  & \\
${\cal C}_6$: & ${\cal C}_5$, & ${\cal C}_4 \sqcup {\cal C}_1$, & ${\cal C}_3 \sqcup {\cal C}_2$, &  & 
\end{tabular}

\bigskip
\noindent 
which implies the following generic hyperplane compositions

\bigskip
\begin{tabular}{llllll}
${\cal C}_N$: & ${\cal C}_{N-1}$, & ${\cal C}_{N-2} \sqcup {\cal C}_1$, & ${\cal C}_{N-3} \sqcup {\cal C}_2$, & \ldots,  & ${\cal C}_{\frac{N}{2}} \sqcup {\cal C}_{\frac{N}{2}-1}$, 
\end{tabular}

\bigskip
\noindent 
or

\bigskip
\begin{tabular}{llllll}
${\cal C}_N$: & ${\cal C}_{N-1}$, & ${\cal C}_{N-2} \sqcup {\cal C}_1$, & ${\cal C}_{N-3} \sqcup {\cal C}_2$, & \ldots, & ${\cal C}_{\lfloor\frac{N}{2}\rfloor} \sqcup {\cal C}_{\lfloor\frac{N}{2}\rfloor}$,  
\end{tabular}

\bigskip
\noindent 
according as $N$ is even or odd, respectively; here, the symbol `$\sqcup$' stands for a disjoint union of two sets.

In the spirit of previous sections, let us also have a closer look at the nature of our generic ${\cal C}_N$. To this end, we first recall the following observations. ${\cal C}_4$, the Desargues configuration, can be viewed as ($4-2=$) two triangles in perspective from a line which are also perspective from a point, that is  ${\cal C}_1$; the line and the point form a geometric hyperplane of ${\cal C}_4$. Next, ${\cal C}_5$, the Cayley-Salmon configuration, admits a view as ($5-2=$) three triangles in perspective from a line where the points of perspectivity of three pairs of them are on a line, {\it aka} ${\cal C}_2$; the two lines form a geometric hyperplane of ${\cal C}_5$. Finally, ${\cal C}_6$, our $(21_5,35_3)$-configuration, can be treated as ($6-2=$) four triangles in perspective from a line where the points of perspectivity of six pairs of them lie on a Pasch configuration, {\it alias} ${\cal C}_3$; the line and the Pasch configuration form a geometric hyperplane of ${\cal C}_6$. Generalizing these observations, we conjecture that for {\it any} $N \geq 4$, ${\cal C}_N$ can be regarded as $N-2$ triangles that are in perspective from a line in such a way that the points of perspectivity of ${N-2 \choose 2}$ pairs of them form the configuration isomorphic to ${\cal C}_{N-3}$, where the latter and the axis of perspectivity form a geometric hyperplane of ${\cal C}_N$. 

\begin{figure}[pth!]
	\centerline{\includegraphics[width=7.0truecm]{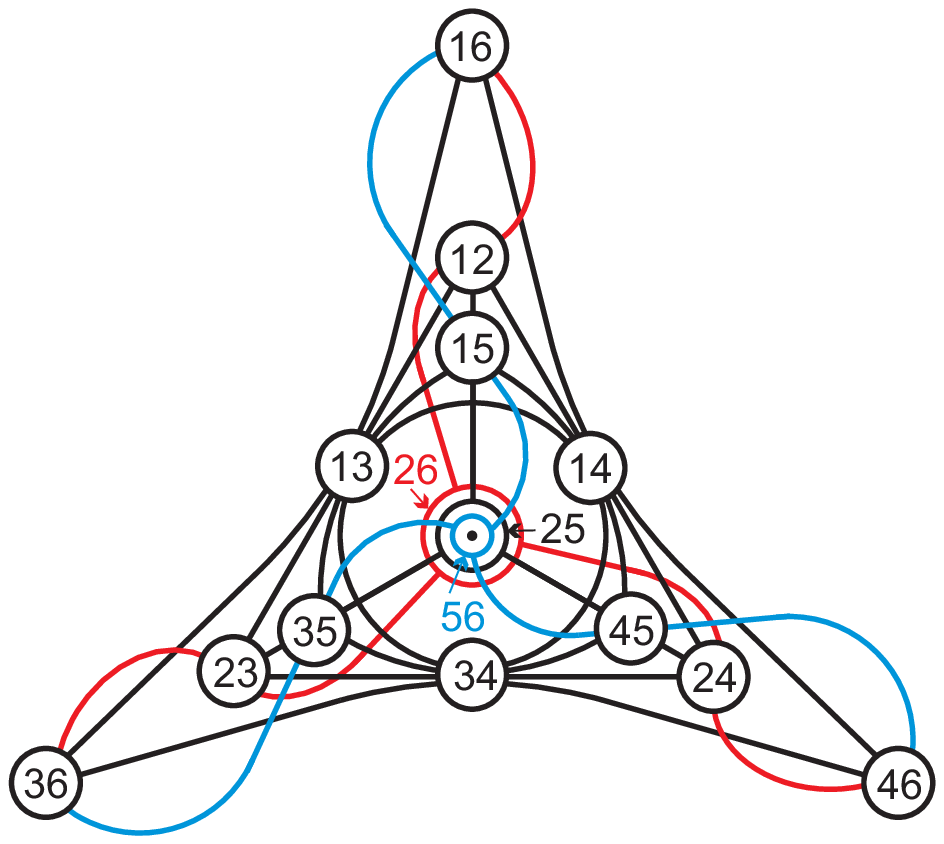}\hspace*{.6cm}\includegraphics[width=6truecm]{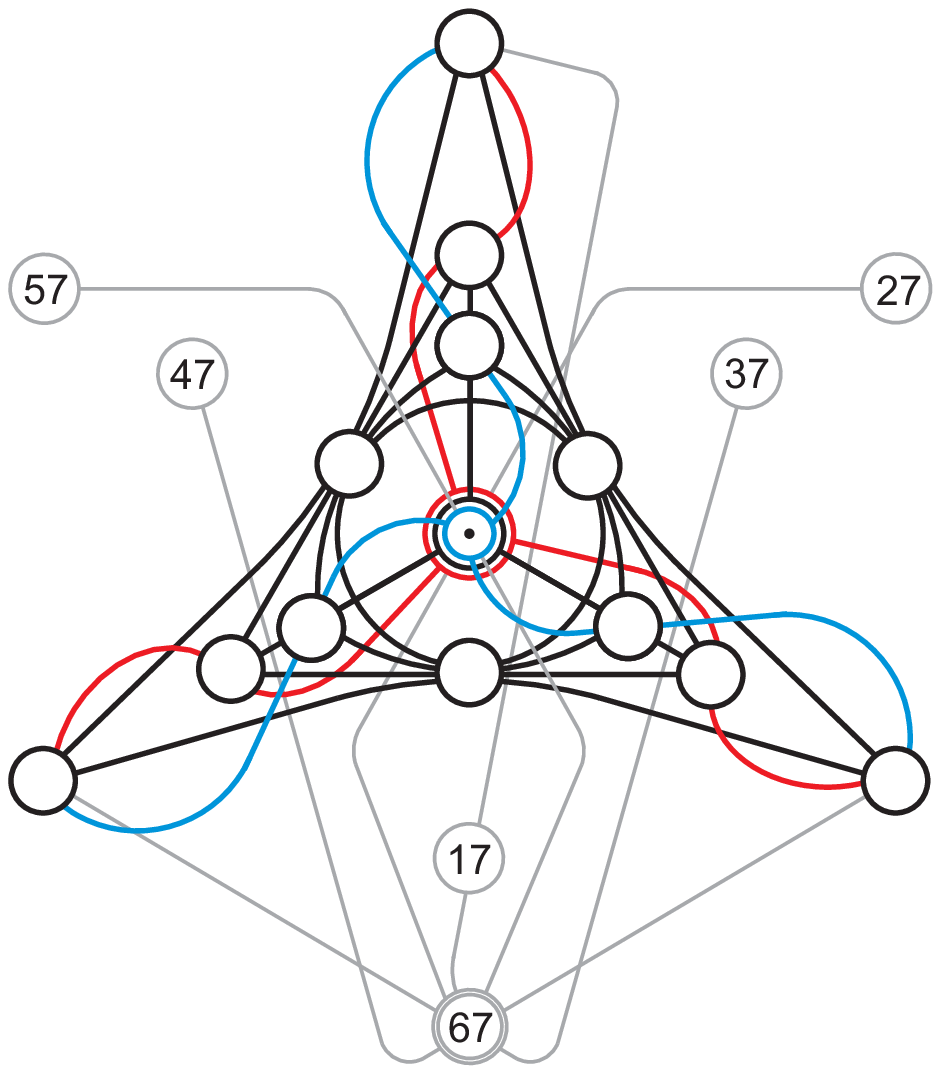}}
	\caption{{\it Left}: -- A diagrammatical proof of the isomorphism between ${\cal C}_5$ and $G_2(6)$.
The points of ${\cal C}_5$ are labeled by pairs of elements from the set $\{1,2,\ldots, 6\}$ in such a way that each line of the configuration is indeed of the form $\{\{a,b\}, \{a,c\}, \{b,c\}\}$, $a \neq b \neq c \neq a$. {\it Right}: -- A pictorial illustration of ${\cal C}_6 \cong G_2(7)$. Here, the labels of six additional points are only depicted, the rest of the labeling being identical to that shown in the left-hand side figure.}
	\label{grass}
	\end{figure}

Next, we invoke the concept of combinatorial Grassmannian (see, e.\,g., \cite{pra1,pra2}). Briefly, a combinatorial Grassmannian $G_k(|X|)$, where $k$ is a positive
integer and $X$ is a finite set, is a point-line incidence structure whose points are $k$-element subsets of $X$ and whose lines are $(k + 1)$-element subsets of $X$, incidence being inclusion. It is known \cite{pra1} that if $|X| =
N+1$ and $k=2$,
$G_2(N+1)$ is a binomial $\left({N+1 \choose 2}_{N-1}, {N+1 \choose
3}_{3}\right)$-configuration; in particular, $G_2(4)$ is the Pasch configuration, 
$G_2(5)$ is the Desargues configuration and $G_2(N+1)$'s with $N \geq 5$ are called {\it
generalized} Desargues configurations. Now, from our detailed examination of the four cases it follows that ${\cal C}_3$, ${\cal C}_4$, ${\cal C}_5$, ${\cal C}_6$, \ldots, ${\cal C}_N$ are endowed with 1, 5, 15, 35, \ldots, ${N+1 \choose 4}$ Pasch configurations. And as ${N+1 \choose 4}$ is also the number of Pasch configurations in $G_2(N+1)$, $N \geq 3$, we are also naturally led to conjecture that ${\cal C}_N \cong G_2(N+1)$.
From what we have found in the previous sections it follows that this property indeed holds for $1 \leq N \leq 6$, being illustrated for 
$N=5$ and $N=6$ in Figure \ref{grass}.

\section{Conclusion}
An intriguing finite-geometrical underpinning of the multiplication tables of Cayley-Dickson algebras $A_N$, $3 \leq N \leq 6$, has been found that admits generalization to any higher-dimensional $A_N$. This started with an observation that the multiplication properties of imaginary units of the algebra $A_N$ are encoded in the structure of the projective space PG$(N-1,2)$. Next, this space was shown to possess a refined structure stemming from particular properties of triples of imaginary units forming its lines. To account for this refinement, we employed the concept of Veldkamp space of point-line incidence structure and found out the latter to be a binomial $\left({N+1 \choose 2}_{N-1}, {N+1 \choose 3}_{3}\right)$-configuration ${\cal C}_N$; in particular, ${\cal C}_3$ (octonions)  was found to be isomorphic to the Pasch $(6_2,4_3)$-configuration, ${\cal C}_4$ (sedenions) to the famous Desargues $(10_3)$-configuration, ${\cal C}_5$ (32-nions) to the Cayley-Salmon $(15_4,20_3)$-configuration found in the well-known Pascal mystic hexagram and ${\cal C}_6$ (64-nions) was shown to be identical with a particular
$(21_5,35_3)$-configuration that can be viewed as four triangles in perspective from a line where the points of perspectivity of six pairs of
them form a Pasch configuration.
These configurations are seen to form a remarkable nested pattern, where ${\cal C}_{N-1}$ is embedded in ${\cal C}_N$ as its geometric hyperplane, that naturally reflects the spirit of the Cayley-Dickson recursive construction of corresponding algebras. 

It is a well-known fact that the only first four algebras $A_N$, $0 \leq N \leq 3$, are `well-behaving' in the sense of being normed, alternative and devoid of zero-divisors --- the facts that are frequently offered as an explanation why a relatively little attention has been paid so far to their higher-dimensional cousins, these latter being even regarded  by some scholars as `pathological.'  It may well be that our finite-geometric, Veldkamp-space-based approach will be able to shed a novel, unexpected light at this issue as it is only starting with $N=4$ when  ${\cal C}_N$ is found to feature a `generalized Desargues property' in the sense that it can be interpreted as $N-2$ triangles that are in perspective from a line in such a way that the points of perspectivity of ${N-2 \choose 2}$ pairs of them form the configuration isomorphic to ${\cal C}_{N-3}$. Or, in  a slightly different form, it is only for $N \geq 4$ when ${\cal C}_N$ contains {\it Desargues} configurations, these occurring as components of its geometric hyperplanes at that.


\section*{Acknowledgments}
This work was partially supported by the VEGA Grant Agency, Project 2/0003/13, as well as by the Austrian Science Fund (Fonds zur F\"orderung der Wissenschaftlichen Forschung (FWF)), Research Project M1564--N27. We would like to express our sincere gratitude to Hans Havlicek (Vienna University of Technology) for a number of critical remarks on the first draft of the paper. 
Our special thanks go to  J\"org Arndt who computed for us multiplication tables of octonions, sedenions, 32-nions and 64-nions. We are also very grateful to Boris Odehnal (University of Applied Arts, Vienna) for creating nice computer versions of several figures.

\end{document}